\documentclass{article}
\usepackage[utf8]{inputenc}

\usepackage[margin=3.2cm]{geometry}

\usepackage{setspace}

\usepackage{graphics}

\usepackage{docmute}

\usepackage[utf8]{inputenc}

\usepackage{mathtools}

\usepackage{amsmath}

\usepackage{amsfonts}

\usepackage{amssymb}

\usepackage{amsthm}
\usepackage{aliascnt}

\usepackage{makeidx}
\makeindex

\usepackage{textcomp}
\usepackage[sb]{libertine}
\usepackage[varqu,varl]{zi4}%
\usepackage[libertine,bigdelims,vvarbb]{newtxmath}
\usepackage[supstfm=libertinesups,supscaled=1.2,raised=-.13em]{superiors}
\usepackage{bm}
\useosf

\usepackage[breaklinks,linktocpage]{hyperref}
\usepackage[hyperpageref]{backref}

\usepackage{url}
\usepackage{breakurl}

\usepackage{tikz-cd}

\usepackage{enumitem}
\usepackage[UKenglish]{babel}
\usepackage[makeroom]{cancel}
\usepackage{pstricks-add}

\definecolor{shadecolor}{rgb}{0.88,0.91,0.95}       
\usepackage{tcolorbox}

\usepackage{dynkin-diagrams}
\usetikzlibrary{backgrounds}

\usepackage{booktabs}

\usepackage{fancyhdr}
\newcommand{\R}{\mathbb{R}}
\newcommand{\C}{\mathbb{C}}

\newcommand{\Ker}{\operatorname{Ker}}

\newcommand{\Ric}{\operatorname{Ric}}

\newcommand{\rank}{\operatorname{rank}}

\newcommand{\CP}{\mathbb{CP}}
\newcommand{\RP}{\mathbb{RP}}

\renewcommand{\|}[1]{\left| \left| #1 \right| \right|}

\newcommand{\<}{\langle}
\renewcommand{\>}{\rangle}

\newcommand{\grad}{\operatorname{grad}}

\newcommand{\sing}{\operatorname{sing}}

\newcommand{\fix}{\operatorname{fix}}

\renewcommand{\d}{\operatorname{d} \!}

\renewcommand{\tilde}[1]{\widetilde{#1}}

\newcommand\Item[1][]{%
  \ifx\relax#1\relax  \item \else \item[#1] \fi
  \abovedisplayskip=0pt\abovedisplayshortskip=0pt~\vspace*{-\baselineskip}}

\usepackage[capitalise]{cleveref}

\numberwithin{equation}{section}

\newtheorem{proposition}[equation]{Proposition}
\crefname{proposition}{Proposition}{Propositions}

\newtheorem{lemma}[equation]{Lemma}
\crefname{lemma}{Lemma}{Lemmas}

\newtheorem{corollary}[equation]{Corollary}
\crefname{corollary}{Corollary}{Corollaries}

\newtheorem*{corollary*}{Corollary}
\crefname{corollary*}{Corollary}{Corollaries}

\crefname{theorem}{Theorem}{Theorems}

\newtheorem{conjecture}[equation]{Conjecture}
\crefname{conjecture}{Conjecture}{Conjectures}

\newtheorem*{theorem*}{Theorem}
\crefname{theorem*}{Theorem}{Theorems}

\crefname{claim}{Claim}{Claims}

\theoremstyle{remark}

\crefname{question}{Question}{Questions}

\crefname{definition}{Definition}{Definitions}

\crefname{example}{Example}{Examples}

\newtheorem{remark}[equation]{Remark}
\crefname{remark}{Remark}{Remarks}

\crefname{assumption}{Assumption}{Assumptions}

\usepackage{pgfplots}
\pgfplotsset{compat=1.15}
\usepackage{mathrsfs}
\usetikzlibrary{arrows}

\usepackage[font={small}]{caption}
\usepackage{authblk}
\usepackage{hyperref}

\title{Harmonic $1$-forms on real loci of Calabi-Yau manifolds}
\author[1]{\small Michael R. Douglas\thanks{\href{mailto:mdouglas@cmsa.fas.harvard.edu}{mdouglas@cmsa.fas.harvard.edu}}}
\author[2,3]{\small Daniel Platt\thanks{\href{mailto:daniel.platt.berlin@gmail.com}{daniel.platt.berlin@gmail.com}}}
\author[4,5]{\small Yidi Qi\thanks{\href{mailto:y.qi@northeastern.edu}{y.qi@northeastern.edu}}}
\author[6]{\small Rodrigo Barbosa\thanks{\href{mailto:rdmbarbosa@gmail.com}{rdmbarbosa@gmail.com}}}

\affil[1]{\small CMSA, Harvard University, Cambridge, MA 02138, USA}
\affil[2]{\small Department of Mathematics, Imperial College London, SW7 2RH London, United Kingdom}
\affil[3]{\small Imperial-X, Imperial College London, W12 0BZ London, United Kingdom}
\affil[4]{\small Department of Physics, Northeastern University, Boston, MA 02115, USA}
\affil[5]{\small NSF Institute for Artificial Intelligence and Fundamental Interactions, Boston, MA, USA}
\affil[6]{\small Revelio Labs, 43 W 23rd St, New York, NY 10010, USA}

\date{\today}

\begin{document}

\maketitle

\abstract{We numerically study whether there exist nowhere vanishing harmonic $1$-forms on the real locus of some carefully constructed examples of Calabi-Yau manifolds, which would then give rise to potentially new examples of $G_2$-manifolds and an explicit description of their metrics.
We do this in two steps:
first, we use a neural network to compute an approximate Calabi-Yau metric on each manifold.
Second, we use another neural network to compute an approximately harmonic $1$-form with respect to the approximate metric, and then inspect the found solution.
On two manifolds existence of a nowhere vanishing harmonic $1$-form can be ruled out using differential geometry.
The real locus of a third manifold is diffeomorphic to $S^1 \times S^2$, and our numerics suggest that when the Calabi-Yau metric is close to a singular limit, then it admits a nowhere vanishing harmonic $1$-form.
We explain how such an approximate solution could potentially be used in a numerically verified proof for the fact that our example manifold must admit a nowhere vanishing harmonic $1$-form.}

\tableofcontents

\def\IC{{\mathbb C}}
\def\IR{\mathbb{R}}
\def\IK{{\mathbb K}}
\def\IZ{\mathbb{Z}}
\def\IP{\mathbb{P}}
\def\E{{\mathbb E}}
\def\Var{\mbox{Var}}
\def\CF{{\mathcal F}}
\def\CH{{\mathcal H}}
\def\CI{{\mathcal I}}
\def\CL{{\mathcal L}}
\def\CM{{\mathcal M}}
\def\CN{{\mathcal N}}
\def\CO{{\mathcal O}}
\def\CP{{\mathcal P}}
\def\CW{{\mathcal W}}
\def\CX{{\mathcal X}}
\def\CY{{\mathcal Y}}
\newcommand{\be}{\begin{equation}}
\newcommand{\ee}{\end{equation}}
\newcommand{\ba}{\begin{array}}
\newcommand{\ea}{\end{array}}
\newcommand{\bea}{\begin{eqnarray}}
\newcommand{\eea}{\end{eqnarray}}
\newcommand{\bigslant}[2]{{\raisebox{.2em}{$#1$}\left/\raisebox{-.2em}{$#2$}\right.}}

\section{Introduction}

Calabi-Yau manifolds are complex manifolds that carry a Ricci-flat Riemannian metric which are immensely important objects in superstring theory as well as pure mathematics.
Thanks to Yau's proof of the Calabi conjecture in \cite{Yau1977,Yau1978} billions of examples of Calabi-Yau manifolds have been constructed using methods in algebraic geometry.
However, not much is known beyond existence of these Ricci-flat metrics.
In fact, it is widely believed that no closed form expression for them exists.

For some applications it would be desirable to have information about the metric.
In Physics, this is needed to compute particle masses and couplings, and very often approximate Calabi-Yau metrics are used to approximately compute these quantities, see e.g. \cite{Ashmore2021,Constantin2024,Butbaia2024}.
In this article, we focus on the following application in differential geometry:
in \cite{Joyce2021} a new construction method for $G_2$-manifolds, which are Ricci-flat Riemannian manifolds in dimension $7$, was explained.
Contrary to Calabi-Yau manifolds, not many constructions for these manifolds are known.
However, they are important objects in M-theory and differential geometry, and new examples would help to generate conjectures tackling open questions in the field, for example which smooth $7$-manifolds can be $G_2$-manifolds.
The construction in \cite{Joyce2021} requires as input a Calabi-Yau manifold $M$ of complex dimension $3$ whose real locus $L \subset M$ (of real dimension $3$) admits a harmonic nowhere vanishing $1$-form $\lambda \in \Omega^1(L)$.
Here, \emph{harmonic} means harmonic with respect to the elusive Calabi-Yau metric.
So far, the nowhere vanishing condition could not be checked in any non-trivial example.
(It could be checked on the complex $3$-dimensional manifold $T^2 \times \text{K3}$, where $T^2$ denotes the torus of real dimension $2$.)
The construction method gives quite an explicit description of the Ricci-flat metric on the $7$-manifold.
This would allow approximately computing some quantities in M-theory on these examples, paralleling the previously described computations on Calabi-Yau manifolds from the physics literature.

In this article, we study the existence of a harmonic nowhere vanishing $1$-form numerically in some examples.
Our motivation for this is twofold:

\begin{enumerate}
    \item 
    There currently exists only one conjectural example for a Calabi-Yau manifold with the required $1$-form in the literature, namely \cite[Section 7.5]{Joyce2021}.
    Our numerical study suggests one new conjectural example.

    \item 
    Our numerical study provides an \emph{approximately harmonic} $1$-form for an \emph{approximate Calabi-Yau metric}.
    Through the advent of greater computational power and specialised software packages it has in recent years become possible to prove existence of genuine solutions to an equation near to a computer-generated approximate solution, see \cite{Gomez-Serrano2019} for an overview of such computer-assisted proofs.
    We note that our approximate solution is likely not of sufficient quality to complete such a computer-assisted proof, but it is possible that our approach could be used to produce suitable approximate solutions in the future.
\end{enumerate}

The article is structured as follows:

\begin{enumerate}
    \item[Section 2.]
    In this section we fix basic notation of Calabi-Yau manifolds and review the need for harmonic nowhere vanishing harmonic $1$-forms based on \cite{Joyce2021}.

    \item[Section 3.]
    In \cref{subsection:construct-nowhere-vanishing-from-approximate} we formalise the second item from above:
    in \cref{corollary:metric-and-1-form-combined-corollary} we show that if one is given a nowhere vanishing harmonic $1$-form for an approximate Calabi-Yau metric, then there exists a nearby nowhere vanishing harmonic $1$-form for the genuine Calabi-Yau metric, provided that the approximate Calabi-Yau metric is close to solving the complex Monge–Ampère equation.
    This uses a standard proof method in geometric analysis.

    In \cref{subsection:candidate-calabi-yaus} we define the Calabi-Yau manifolds on which we run our experiments:
    the only harmonic $1$-form for \texttt{Fermat} is constant zero, for \texttt{Quintic} any harmonic $1$-form is known to have a zero thanks to some results about the geometry of $3$-manifolds.
    We still run our numerics on them as a control.
    The real locus of \texttt{CICY1} and \texttt{CICY2} is diffeomorphic to $S^1 \times S^2$, so might admit a nowhere vanishing harmonic $1$-form.

    \item[Section 4.]
    In this section we explain our numerical algorithms:
    we often numerically evaluate integrals, and this is done using Monte-Carlo integration, which is explained in \cref{subsection:sampling-algorithms}.
    In \cref{section:bihomogeneous-nn} we explain how our approximate Calabi-Yau metrics are computed:
    we review the approach from \cite{Douglas2022} and explain our extension to higher codimension.
    Following this, we explain in \cref{section:numerical-forms} how we compute approximately harmonic $1$-forms with respect to the approximate Calabi-Yau metric on the real locus of a Calabi-Yau manifold.

    \item[Section 5.]
    Finally, in this section we present our numerical results.
    \Cref{subsection:results-approximate-cy-metrics} discusses properties of the approximate Calabi-Yau metrics on our example varieties.
    Our data supports the folklore that the approximate metrics are of better quality in regions of the manifold that are far away from a smoothing of singularities.
    We also experimentally observe the formation of long necks that for hypersurfaces is guaranteed by \cite[Theorem 1.2]{Sun2019}, but in higher codimension is not known.

    Most importantly, in \cref{section:results-1-forms} we discuss the approximate $1$-forms.
    \texttt{Fermat} and \texttt{Quintic} cannot have nowhere vanishing harmonic $1$-forms which is in agreement with our experimental results.
    The real locus of \texttt{CICY2} is diffeomorphic to $S^1 \times S^2$, and the approximate harmonic $1$-form is approximately constant in the $S^1$-direction and zero in the $S^2$-direction.
    This can be seen from quantitative analysis of our approximate solution, but is even clearly visible when plotting the $1$-form on $1$-dimensional slices of the real locus.    
\end{enumerate}

The code pertaining to our experiments is available at \url{https://github.com/yidiq7/MLHarmonic-1-form}.

\textbf{Acknowledgments.}
The second author wishes to thank Oliver Gäfvert, Michel van Garrel, Andre Lukas, Lorenzo la Porta, Simon Salamon, Calum Spicer, and Mehdi Yazdi for helpful conversations. 
The third author would like to thank Simon Salamon, Challenger Mishra and Justin Tan for their hospitality during his visit at King's College London and University of Cambridge, where part of this work was carried out. He also thanks James Halverson, Edward Hirst, Dami\'an Mayorga Pe\~na, and his supervisor Fabian Ruehle for useful discussion. The third author acknowledges support from the National Science Foundation under Cooperative Agreement PHY-2019786 (the NSF AI Institute for Artificial Intelligence and Fundamental Interactions).
The fourth author would like to thank Simon Donaldson for helpful conversations, and to acknowledge support from the Simons Center for Geometry and Physics.
The authors thank the anonymous referee for their helpful comments.

\section{Background}
\subsection{Complete intersection Calabi-Yau manifolds}

We begin by surveying some basic facts about Calabi-Yau manifolds, as can be found in \cite{Yau2009,Sun2022}.
A Kähler manifold $(M^{2m},J,\omega,g)$, where $J$ is a complex structure, and $\omega$ and $g$ are a compatible symplectic form and Riemannian metric respectively, is called \emph{Calabi-Yau}, if there exists a holomorphic $m$-form $\Omega \in \Omega^{n,0}(M)$ that is nowhere vanishing.
For simplicity, we assume that $\int_M \omega^m=\int_M \Omega \wedge \overline{\Omega}$.
In \cite{Yau1977,Yau1978}, Yau proved the Calabi conjecture, showing in particular that on a Calabi-Yau manifold, there exists a Kähler form $\omega_{\text{CY}} \in \Omega^2(M)$ in the same cohomology class as $\omega$, such that 
\begin{align}
\label{equation:cy-kaehler-and-holo-volume-form-comparison}
\omega_{\text{CY}}^m=\Omega \wedge \overline{\Omega}.
\end{align}
The corresponding metric $g_{\text{CY}}$ is then Ricci-flat.

This reduces the task of constructing Ricci-flat metrics to finding complex Kähler manifolds with holomorphic volume forms.
Using algebraic geometry techniques, billions of examples satisfying this criterion have been constructed.
A large class of 7890 examples comes from \emph{Complete Intersection Calabi-Yau} manifolds (CICYs), as worked out in \cite{Candelas1987}.
(A complete list of examples can be accessed at \url{https://www-thphys.physics.ox.ac.uk/projects/CalabiYau/cicylist/}.)

\subsection{$G_2$-orbifolds and their resolutions}

The first examples of compact manifolds with holonomy $G_2$ were constructed in \cite{Joyce1996}.
These manifolds were obtained as the resolution of orbifolds $T^7/\Gamma$, where $\Gamma$ is a group of involutions acting on $T^7$.

In \cite{Joyce2021}, this method was generalised to orbifolds of the form $M/\Gamma$, where $M$ is any $G_2$-manifold and $\Gamma$ is a group of involutions preserving the $G_2$-structure, acting on $M$.
Choosing $M=T^7$ recovers Joyce's original construction.

We denote $L:=\fix(\Gamma)$.
By \cite[Prop. 10.8.1]{Joyce2000}, the submanifold $L$ is a three-dimensional associative submanifold.
The $G_2$-manifold $N$ is constructed by gluing an Eguchi-Hanson space into $M/\Gamma$ over every point in $L$.
The construction requires a harmonic $1$-form $\lambda \in \Omega^1(L)$ which defines the size of the glued in Eguchi-Hanson space.
Crucially, the construction only works if $\lambda$ is nowhere $0$.
It is an open problem to define a suitable $G_2$-structure around a zero of $\lambda$, see \cite[Section 8, point (vi)]{Joyce2021}.

Currently, the only known examples that produce manifolds with holonomy equal to $G_2$ are:
\begin{enumerate}
    \item 
    $M=T^7$,

    \item 
    $M=T^3 \times X$, where $X$ is a K3 surface, and $L=S^1 \times S^2 \cup S^1 \times C$ for a holomorphic curve $C \subset X$.
    Here, the metric on $L$ is approximately the product metric, and a harmonic $1$-form that is approximately the parallel $1$-form in the $S^1$-direction is known to exist, see \cite[Section 7.3]{Joyce2021}.
\end{enumerate}
It is believed that many more examples of the form
\[
    M=S^1 \times Y
\]
exist, where $(Y,J,\omega,g)$ is a Calabi-Yau manifold of complex dimension $3$ with an anti-holomorphic involution $\sigma: Y \rightarrow Y$.
As before, denote by $g_{\text{CY}}$ the Calabi-Yau metric in the same Kähler class as $g$.
If we define $\hat{\sigma}:M \rightarrow M$, $(t,y) \mapsto (-t, \sigma(y))$, and if $L':=\fix(\sigma)$ admits a nowhere vanishing $1$-form that is harmonic with respect to $g_{\text{CY}}$, then $M/\< \hat{\sigma}\>$ can be resolved and one obtains a manifold with holonomy $G_2$ (if Y had holonomy $\text{SU}(3)$), see \cite[Section 7.5]{Joyce2021}.
However, no such example has been proven to exist, because it is hard to check if $\lambda$ is nowhere vanishing, because $g_{\text{CY}}$ is not known explicitly.

\section{Calabi-Yau manifolds and harmonic $1$-forms}
\subsection{Constructing nowhere vanishing harmonic $1$-forms from approximate harmonic $1$-forms}
\label{subsection:construct-nowhere-vanishing-from-approximate}

Let $(M,J,\omega,g)$ be a Calabi-Yau manifold with Calabi-Yau metric $g_{\text{CY}}$ and $L=\fix(\sigma)$ for an involution $\sigma: M \rightarrow M$ that is anti-holomorphic.
It is usually easy to compute the first Betti number of $L$, and if this is bigger than $0$, then it is known to admit a harmonic $1$-form.
However, it has been lamented above that it is difficult to check if such a $1$-form is \emph{nowhere vanishing}.

One potential avenue to checking this is the following three-step process:
\begin{enumerate}
    \item 
    Compute an approximation $g$ of the unknown Calabi-Yau metric $g_{\text{CY}}$.
    How to do this is explained in \cref{section:bihomogeneous-nn}.

    \item 
    Without explicitly computing $g_{\text{CY}}$, one can estimate the difference $\|{g-g_{\text{CY}}}_{C^{1,\alpha}}$ by computing a certain loss for $g$.
    This is explained in \cref{corollary:g-approx-minus-cy-estimate}.

    \item 
    If the difference $\|{g-g_{\text{CY}}}_{C^{1,\alpha}}$ is sufficiently small, and there exists a harmonic $1$-form for $g$ that is sufficiently bounded from below, then that implies that also $g_{\text{CY}}$ must admit a harmonic $1$-form that is nowhere vanishing.
    This is a consequence of \cref{proposition:perturbation-of-harmonic-1-forms}.
\end{enumerate}

Thus, one can establish existence of a nowhere vanishing harmonic $1$-form with respect to $g_{\text{CY}}$ without computing $g_{\text{CY}}$, just by doing explicit calculations for the approximate metric.
This is called a \emph{numerically verified proof}.
An example in geometry is \cite{Reiterer2019}, other examples are given in \cite[Example 1.2]{Gomez-Serrano2019}.

We point out that the main difficulty in this type of numerically verified proof lies in explicitly bounding the inverse of linearisation of the equation to be solved.
We do not give such an explicit bound here and we provide no numerically verified proof.
This section serves as motivation for our numerical experiments.

The following proposition shows how one obtains a $C^{3,\alpha}$-estimate for the solution of the complex Monge–Ampère equation \cref{equation:cx-monge-ampere}.

\begin{proposition}
    \label{proposition:g-approx-minus-cy-estimate}
    Let $(M,J)$ be a compact, complex manifold of dimension $m$, and $g$ be a Kähler metric with Kähler form $\omega$.
    Then, for all $\epsilon > 0$ there exists $\delta > 0$ (depending on $g$) such that the following is true:

    if $(e^f-1) \in C^{1,\alpha}(M)$ has mean zero and satisfies $\|{e^f-1}_{C^{1,\alpha}} \leq \delta$, then the equation 
    \begin{align}
        \label{equation:cx-monge-ampere}
        (\omega+\d \d^c \phi)^m=e^f \omega^m
    \end{align} 
    has a solution satisfying $\|{\phi}_{C^{3,\alpha/2}} \leq \epsilon$.
\end{proposition}

\begin{proof}
    The linearisation of \cref{equation:cx-monge-ampere} is $\phi \mapsto * \Delta \phi$.
    Assume the operator norm of the operator $\Delta^{-1}: C^{1,\alpha}_0(M) \rightarrow C^{3,\alpha}_0(M)$, mapping mean-zero functions to mean-zero functions, satisfies $\|{\Delta^{-1}} \leq C$.
    Define a sequence $\phi_k$ as follows:
    let $\phi_0=0$, and let $\phi_{k+1}$ be the unique mean-zero solution to
    \begin{align}
    \label{equation:complex-monge-ampere-reordered}
        * \Delta \phi_{k+1}
        =
        -(\d \d^c \phi_k)^2 \wedge \omega^{m-2}
        -
        \dots
        -(\d \d^c \phi_k)^{m-1} \wedge \omega
        -(\d \d^c \phi_k)^m
        +
        (e^f-1)\omega^m.
    \end{align}
    (Note that this equation is essentially just \cref{equation:cx-monge-ampere} reordered so that the linear term is isolated on the left hand side.)
    For simplicity, write $K$ for a number satisfying $|\omega^m| \leq K$ pointwise everywhere.
    Let $\delta$ s.t. $4\delta^2 C^2K^2m<\delta$.
    Then we can prove via induction that $\|{\phi_k}_{C^{3,\alpha}} \leq 2\delta CK$ for all $k \geq 0$, because:
    \begin{align*}
        &\|{\phi_{k+1}}_{C^{3,\alpha}}
        \\
        &\leq
        C\|{\Delta \phi_{k+1}}_{C^{1,\alpha}}
        \\
        &
        \leq
        C
        \left(
        \|{\d \d^c \phi_k)^2 \wedge \omega^{m-2}}_{C^{1,\alpha}}
        +
        \dots
        +
        \|{(\d \d^c \phi_k)^{m-1} \wedge \omega}_{C^{1,\alpha}}
        +
        \|{(\d \d^c \phi_k)^m}_{C^{1,\alpha}}
        +
        \|{(e^f-1)\omega^m}_{C^{1,\alpha}}
        \right)
        \\
        &\leq
        CK
        \left(
        (2\delta CK)^2
        +
        \dots
        +
        (2\delta CK)^{2m-2}
        +
        (2\delta CK)^{2m}
        +
        \delta
        \right)
        \\
        &\leq
        CK
        \left(
        m(2\delta CK)^2
        +
        \delta
        \right)
        \\
        &\leq
        2\delta CK.
    \end{align*}
    By the Arzelà-Ascoli theorem, we therefore have that $\phi_k$ has a subsequence converging in $C^{3,\alpha/2}$.
    By taking the limit in \cref{equation:complex-monge-ampere-reordered}, we see that this limit solves \cref{equation:cx-monge-ampere}.
    And by taking $\delta \leq \frac{\epsilon}{2CK}$ we have that $\|{\phi_k}_{C^{3,\alpha}} \leq 2\delta CK \leq \epsilon$, proving the claim.
\end{proof}

\begin{remark}
    Note that \cref{proposition:g-approx-minus-cy-estimate} does \emph{not} give an alternative proof of the Calabi conjecture.
    In the proposition, existence of $\phi$ is only established if the initial Kähler form $\omega$ is already very close to solving the complex Monge–Ampère equation \cref{equation:cx-monge-ampere}.
\end{remark}

An immediate corollary is the following:

\begin{corollary}
    \label{corollary:g-approx-minus-cy-estimate}
    Let $(M,J,\omega,g)$ be a Calabi-Yau manifold with holomorphic volume form $\Omega \in \Omega^{n,0}(M)$ satisfying $\int_M \omega^m=\int_M \Omega \wedge \overline{\Omega}$ and Calabi-Yau metric $g_{CY}$.
    Then, for all $\epsilon >0$ there exists $\delta >0$ (depending on $g$) such that the following is true:

    Let $f \in C^\infty(M)$ be defined by $e^f \omega^m=\Omega \wedge \overline{\Omega}$.
    If $\|{e^f-1}_{C^{1,\alpha}} \leq \delta$, then $\|{g-g_{\text{CY}}}_{C^{1,\alpha/2}} \leq \epsilon$.
\end{corollary}

This corollary can, in theory provide an estimate for the difference $g-g_{\text{CY}}$.
In practice, a difficulty is to obtain an explicit $\delta$ for a given $\epsilon$.

Roughly speaking, if one is given a nowhere vanishing $1$-form $\lambda \in \Omega^1(L)$ with respect to $g|_L$, and the difference $g-g_{\text{CY}}$ is very small, then a small perturbation of $\lambda$ will be harmonic with respect to $g_{\text{CY}} | _L$.
This is just a fact about Riemannian geometry and has nothing to do with Calabi-Yau metrics, and the following proposition states this for general Riemannian metrics $g$ and $h$:

\begin{proposition}
    \label{proposition:perturbation-of-harmonic-1-forms}
    Let $L$ be a smooth, compact manifold and let $h$ be a Riemannian metric on $L$.
    For all $\epsilon>0$ there exists $\delta>0$ such that the following is true:
    
    Let $g$ be Riemannian metrics on $L$.
    Let $\lambda \in \Omega^1(L)$ be harmonic with respect to $g$ and
    \[
        \|{\lambda}_{L^2,g}=1,
        \text{ and }
        \min_{x \in L} |\lambda|_g(x)
        = \epsilon.
    \]
    There exists $\eta \in \Omega^0(L)$, unique up to constants, such that $\tilde{\lambda}=\lambda+\d \eta$ is harmonic with respect to $h$.
    If $\|{g-h}_{C^1,g} \leq \delta$, then $\tilde{\lambda}$ has no zeros.
\end{proposition}

\begin{proof}
    Assume $\min_{x \in L} |\lambda|_g(x) = \epsilon$ and $\|{g-h}_{C^1,g} \leq \delta$ with $\delta$ to be determined later.
    The operators $\Delta$ and $\d^*$ depend on the metric, and we denote the operators corresponding to $g$ and $h$ by $\Delta_g$, $\d^*_g$ and $\Delta_h$, $\d^*_h$ respectively.

    We have $\Delta_g \lambda=0$, which implies $\d \lambda=0$ and $\d^*_g \lambda =0$ because $L$ is compact.
    Therefore, for all $p \in [1, \infty)$,
    \begin{align}
        \label{equation:d-star-lambda-small}
        \|{
        \d^*_h \lambda
        }_{L^p,h}
        =
        \|{
        (\d^*_h-d^*_g) \lambda
        }_{L^p,h}
        \leq
        \|{
        h-g
        }_{C^1,h}
        \cdot
        \|{
        \lambda
        }_{L^p,h}
        \lesssim
        \delta.
    \end{align}
    Here, in the last step we used that $\|{g-h}_{C^1,g} \leq \delta$ implies $\|{g-h}_{C^1,h} \lesssim \delta$ if $\delta<\frac{1}{2}$, which we assumed without loss of generality.
    We also used 
    $  \|{
        \lambda
        }_{L^p,h} \lesssim 
        \|{
        \lambda
        }_{L^2,h} \lesssim 
        \|{
        \lambda
        }_{L^2,g} =1.$
    In this, the first inequality is true because $\Ker(\Delta_g)$ is a finite-dimensional vector space, so on it all norms are equivalent;
    the second inequality follows from $\|{g-h}_{C^1,h} \lesssim \delta$ with the assumption that $\delta<\frac{1}{2}$.
    Now let $\eta \in \Omega^0(L)$ be the unique function with mean zero satisfying the equation
    \[
        \Delta_h \eta
        =
        -\d^*_h \lambda.
    \]
    Then $\lambda+\d \eta \in \Omega^1(L)$ is the unique $\Delta_h$-harmonic representative in the de Rham cohomology class of $\lambda$, as the following calculation shows:
    \[
        \Delta_h(\lambda+\d \eta)
        =
        \d \d^*_h \lambda+ 
        \underbrace{\d^*_h\d \lambda}_{=0}+
        \underbrace{\d^*_h \d \d \eta }_{=0}+
        \d \d^*_h \d \eta
        =
        \d \d^*_h(\lambda)+\d \, (\Delta_h \eta)
        =
        \d \left(
        \d^*_h \lambda+\Delta_h \eta)
        \right)
        =0.
    \]
    In the second step we used that $\d \d^*_h \d \eta=\d \d^*_h \d \eta+\d \d \d^*_h \eta=\d \Delta_h \eta$.
    We now show that $\d \eta$ is small, and from this it will follow that $\lambda+\d \eta$ is nowhere vanishing.
    Choose $p \in [1, \infty)$ and $\alpha \in (0,1)$ so that $L^p_2 \hookrightarrow C^{1,\alpha}$ is continuously embedded.
    As $L$ is $3$-dimensional, this is for example satisfied for $\alpha=\frac{1}{2}$ and $p=6$.
    Then:
    \begin{align*}
        \|{\d \eta}_{C^0,h}
        \leq
        \|{\eta}_{C^{1,\alpha},h}
        \lesssim
        \|{\eta}_{L^p_2,h}
        \lesssim
        \|{\Delta_h \eta}_{L^p,h}
        =
        \|{\d^*_h \lambda}_{L^p,h}
        \lesssim
        \delta.        
    \end{align*}
    Here, we used the embedding $L^p_2 \hookrightarrow C^{1,\alpha}$ in the second step, we used elliptic regularity of the Laplacian applied to functions orthogonal to the kernel of the Laplacian in the third step, used the definition of $\eta$ in the fourth step, and used \cref{equation:d-star-lambda-small} in the last step.
    This implies
    \[
        \|{\d \eta}_{C^0,g}
        \lesssim
        \|{\d \eta}_{C^0,h}
        \lesssim
        \delta,
    \]
    where we again assumed $\delta<\frac{1}{2}$.
    Spelled out, this means that there exists a constant $C$ depending only on $h$ and the assumption that $\delta<\frac{1}{2}$ but not on $g$ or $\delta$, such that $\|{\d \eta}_{C^0,g} \leq C \delta$.
    Now, if $\delta<\epsilon/C$, then
    \[
        \min_{x \in L}
        |\lambda+\d \eta|_g(x)
        \geq
        \min_{x \in L}
        |\lambda|_g(x)
        -
        \max_{x \in L}
        |\d \eta|_g(x)
        >
        \epsilon-\epsilon
        =0.       
    \]
    This means that $\tilde{\lambda}:=\lambda+\d \eta$ has no zeros.
\end{proof}

In our application, the role of $g$ will be taken by the approximate Calabi-Yau metric, and the role of $h$ will be taken by the Calabi-Yau metric $g_{\text{CY}}$.
So, by combining \cref{corollary:g-approx-minus-cy-estimate,proposition:perturbation-of-harmonic-1-forms}, we obtain:

\begin{corollary}
    \label{corollary:metric-and-1-form-combined-corollary}
    Let $(M,J,\omega,g)$ be a Calabi-Yau manifold with holomorphic volume form $\Omega \in \Omega^{n,0}(M)$ satisfying $\int_M \omega^m=\int_M \Omega \wedge \overline{\Omega}$ and Calabi-Yau metric $g_{CY}$.
    Let $L \subset M$ be the real locus of $M$.
    Then, for all $\epsilon >0$ there exists $\delta >0$ (depending on $g$) such that the following is true:

    Let $f \in C^\infty(M)$ be defined by $e^f \omega^m=\Omega \wedge \overline{\Omega}$.
    If $\|{e^f-1}_{C^{1,\alpha}} \leq \delta$ and $\lambda \in \Omega^1(L)$ is harmonic with respect to $g$ and
    \[
        \|{\lambda}_{L^2,g|_L}=1,
        \text{ and }
        \min_{x \in L} |\lambda|_{g|_L}(x)
        = \epsilon
    \]
    then there exists $\eta \in \Omega^0(L)$, unique up to constants, such that $\tilde{\lambda}=\lambda+\d \eta$ is harmonic with respect to $g_{CY}|_L$ and has no zeros.
\end{corollary}

\subsection{Candidate Calabi-Yau $3$-folds}
\label{subsection:candidate-calabi-yaus}

A necessary condition for $L$ to admit a nowhere vanishing harmonic $1$-form is that it admit a harmonic $1$-form at all, i.e. its first Betti number should satisfy $b_1(L) \neq 0$.
One example is the tetra-quadric from \cite[Example 7.6]{Joyce2021} and \cite{Buchbinder2014}.
In this section, we present one non-example and two other examples.

\subsubsection{The Fermat quintic \texttt{Fermat} with $b_1(L)=0$}
\label{section:fermat}

Let
\[
    \text{\texttt{Fermat}}
    =
    \{ (x_0:\dots:x_4) \in \mathbb{CP}^4:
    x_0^5+\dots+x_4^5=0
    \}
    \subset \mathbb{CP}^4.
\]
This is a simple example of a Calabi-Yau manifold.
The complex conjugation $\sigma: \mathbb{CP}^4 \rightarrow \mathbb{CP}^4$, $(x_0:\dots:x_4) \mapsto (\overline{x_0}:\dots:\overline{x_4})$ restricts to \texttt{Fermat}, and we can study its real locus $L:= \fix(\sigma) \cap \text{\texttt{Fermat}}$.

\begin{lemma}
    The real locus $L$ of \texttt{Fermat} is diffeomorphic to $\RP^3$.
\end{lemma}

\begin{proof}
    The map
    \begin{align*}
        \RP^3 & \rightarrow L
        \\
        (x_0:x_1:x_2:x_3)
        &\mapsto
        \left(
        x_0:x_1:x_2:x_3:
        -\sqrt[5]{x_0^5+x_1^5+x_2^5+x_3^5}
        \right)
    \end{align*}
    is a diffeomorphism.
\end{proof}

Thus, $b_1(L)=0$, and there cannot exist a nowhere vanishing harmonic $1$-form on $L$.
We still run our experiments on \texttt{Fermat} alongside other candidate Calabi-Yau manifolds to make sure our numerical approximations reproduce this fact in order to validate our experimental setup.

\subsubsection{The quintic \texttt{Quintic} with $b_1(L) = 1$}
\label{section:quintic}

We aim to define an explicit quintic that is smooth in $\mathbb{CP}^4$ and whose real locus in $\mathbb{RP}^4$ has $b_1 \neq 0$.
The construction proceeds in three steps:

\begin{enumerate}
\item
Constructing a singular cubic three-fold $D^* \subset \mathbb{RP}^4$ with a single ordinary double point $a$ as its singularity.

\item
Define a smoothing $D$ of $D^*$ with $b_1(D)=1$.

\item
Given $D=Z_{\R}(f)$ we obtain the singular quintic $Q^*=Z_{\R}((x_0^2+\dots+x_4^2)f)$ which is bijective to $D$.
A suitable smoothing $Q$ of $Q^*$ will then still satisfy $b_1(Q)=1$.
\end{enumerate}

We now carry out the three steps in detail.

\paragraph{Constructing a singular cubic with the right type of singularity}
It was shown in the work \cite{Krasnov2006} that the space of smooth cubic surfaces in $\mathbb{RP}^4$ consists of $8$ chambers, separated by $14$ walls.
The exact arrangement of chambers and walls is given in \cite[p.843]{Krasnov2009}.
All walls are sets of singular cubics with a single ordinary double point, and the walls can be distinguished using information coming from the ordinary double point.

We now define a singular cubic $D^*$ in one of these walls, namely $\mathscr{D}^3_1(1)$ in the notation from the article \cite{Krasnov2009}.
The notation $\mathscr{D}^3_1(1)$ means that the bilinear form defining the ordinary double point of $D^*$ has signature $(3,1)$, and that the tautological curve consisting of all lines in $\mathbb{P}(T_a \mathbb{RP}^4)$ that lie completely in $D^*$ has exactly one connected component.

Let
\[
f^*
:=
x_0(x_1^2+x_2^2+x_3^2-x_4^2)-
\left(
x_1^3+x_2^3+x_3^3-
\frac{1}{2}x_4^3
\right)
\]
One checks that $D^*:=Z_{\R}(f^*)$ is a singular variety whose unique singular point is $a=[1:0:0:0:0] \in Z_{\R}(f^*) \subset \mathbb{RP}^4$.
Near $a$, the cubic monomials in the definition of $f^*$ are very small compared to the rest, so $D^*$ locally looks like the ordinary double point $Z_{\R}(x_1^2+x_2^2+x_3^2-x_4^2)$.
The signature of the bilinear form defining this quadric is $(3,1)$, so $D^* \in \mathscr{D}^3_1$, in the notation of \cite{Krasnov2009}.
It remains to check that $D^* \in \mathscr{D}^3_1(1)$.
If we canonically identify $\mathbb{P}(T_a \mathbb{RP}^4)$ with the set of lines in $\mathbb{RP}^4$, then the definition of the tautological curve $C \subset \mathbb{P}(T_a \mathbb{RP}^4)$ is the set of all lines tangent to $D^*$ that lie completely in $D^*$.
We will show that $C$ has exactly one connected component.
The singular point $a$ is contained in the affine chart $\{ x_0 \neq 0 \} \subset \mathbb{RP}^4$, so we can canonically identify $T_a \mathbb{RP}^4= \{(x_1,x_2,x_3,x_4) \in \R^4\} = \R^4$.
The curve is then the following subset in the projectivisation of this space:
\[
C
=
Z_{\R}
\left(
x_1^2+x_2^2+x_3^2-x_4^2, \quad
f^*
\right)
=
Z_{\R}
\left(
x_1^2+x_2^2+x_3^2-x_4^2, \quad
x_1^3+x_2^3+x_3^3-\frac{1}{2}x_4^3
\right)
\subset
\mathbb{RP}^3,
\]
where for the second equality we plugged the equation $x_1^2+x_2^2+x_3^2-x_4^2=0$ into $f^*$.
Now, $C$ is easily seen to be diffeomorphic to $S^1$, for example by plotting the two surfaces $Z_{\R}(x_1^2+x_2^2+x_3^2-x_4^2), Z_{\R}(x_1^3+x_2^3+x_3^3-\frac{1}{2}x_4^3) \subset \mathbb{R}^3$.
Thus, $D^* \in \mathscr{D}^3_1(1)$.

\paragraph{Finding a smoothing with the right topology}
By \cite[p.843]{Krasnov2009}, when smoothing $D^*$, the resulting smooth cubic $D$ has one of two topological types:
\begin{itemize}
    \item 
    either $D \in \mathscr{B}(0)$, i.e. $D \subset \RP^4$ is diffeomorphic to $\RP^3$ and therefore $b_1(D)=0$,

    \item 
    or $D \in \mathscr{B}(1)_{\text{II}}$, i.e. $D \subset \RP^4$ is diffeomorphic to $\R P^3 \# (S^1 \times S^2)$ and therefore $b_1(D)=1$.
\end{itemize}

In $\mathbb{R}^3$, the smoothing $Z_{\R} \left( x_1^2 + x_2^2 - x_2^2 - \frac{1}{2} \right)$ of the ordinary double point $Z_{\R} \left( x_1^2 + x_2^2 - x_2^2\right)$ has one additional non collapsible circle.
Likewise, if one defines $f := f^* - \frac{1}{4} x_0^3$, then $D := Z_{\R}(f) \subset \mathbb{RP}^4$ has an additional non collapsible circle, and therefore must be homeomorphic to $\R P^3 \# (S^1 \times S^2)$.

Note that in this case it is known that every harmonic $1$-form on $D$ must have at least one zero because of the following proposition:

\begin{proposition}
    \label{proposition:must-have-zeros-on-RP3-with-handle}
    Every closed $1$-form on $\R P^3 \# (S^1 \times S^2)$ must have at least one zero.
    If its zeros are non-degenerate (i.e. near each zero it is locally of the form $\d f$ for a function $f$ with non-degenerate Hessian), then it must have an even number of zeros.
\end{proposition}

\begin{proof}
    The manifold $\R P^3 \# (S^1 \times S^2)$ is \emph{reducible} in the sense of \cite[p.18]{Jaco1980}.
    By \cite[VI.7. Lemma]{Jaco1980}, $\R P^3 \# (S^1 \times S^2)$ is not Seifert fibered, in particular it is not a fibration over $S^1$.
    By \cite{Tischler1970}, $\R P^3 \# (S^1 \times S^2)$ does not admit a nowhere vanishing closed $1$-form.
    In particular, a harmonic $1$-form must have at least one zero.
    As $\R P^3 \# (S^1 \times S^2)$ is orientable, its Euler characteristic $\chi(\R P^3 \# (S^1 \times S^2))$ is zero.
    By the Poincaré–Hopf theorem (see e.g. \cite[p.35]{Milnor1965}), this implies that a harmonic $1$-form with non-degenerate zeros must have an even number of zeros.
\end{proof}

\paragraph{Using the cubic to define a quintic}
Define the quintic polynomial $g^*:=(x_0^2+\dots+x_4^2)f$.
Notice first that the factor $(x_0^2+\dots+x_4^2)$ is nowhere vanishing on $\mathbb{RP}^4$.
Thus, $Z_{\R}(f)=Z_{\R}(g^*) \subset \mathbb{RP}^4$.
The real locus of the complex variety $X^* := Z_{\C}(g^*)$ is smooth, but for our application we require an algebraic variety that is smooth as a \emph{complex} variety.
The variety $X^* \subset \mathbb{CP}^4$ is reducible, so in particular singular, but a small generic perturbation of it is smooth.
As $Z_{\R}(g^*) \subset \RP^4$ is smooth, any small perturbation of it does not change its diffeomorphism type.
Altogether, we have shown the following:

\begin{proposition}
\label{proposition:quintic-construction}
There exists $\epsilon_0 > 0$ such that the quintic polynomial
\[
Q=(x_0^2+\dots+x_4^2)
\left(
x_0(x_1^2+x_2^2+x_3^2-x_4^2)-(x_1^3+x_2^3+x_3^3-
\frac{1}{2}x_4^3
)
-\frac{1}{4}x_0^3
\right)
+\epsilon P,
\]
where $\epsilon \in \R$ such that $0<\epsilon < \epsilon_0$ and $P \in \R[x_0,x_1,x_2,x_3,x_4]$ is a generic homogeneous polynomial of degree $5$, defines a smooth algebraic variety in $\mathbb{CP}^4$, and its real locus has first Betti number $b^1=1$.
\end{proposition}

For concreteness, here is how one gets an explicit estimate for $\epsilon_0$:

\begin{proposition}
    \label{proposition:conrete-small-deformation-bound}
    Fix homogeneous degree $5$ polynomials $g^*,P$ such that $Z_{\R}(g^*) \subset \mathbb{RP}^4$ is smooth.
    Let $U \subset S^4$ be a neighbourhood of $Z_{\R}(g^*) \cap S^4 \subset \mathbb{R}^5$.
    Let $\tilde{g}:=g^*|_{S^4}, \tilde{P}:=P|_{S^4} \in C^{\infty}(S^4)$.
    Define
    \[
        m := 
        \min_{x \in U}
        | \nabla \tilde{g} |,
        \quad
        k :=
        \min_{x \in S^4 \setminus U}
        | \tilde{g} |.
    \]
    If $\|{\tilde{P}}_{C^0} < k$ and $\|{\nabla \tilde{P}}_{C^0} < m$, then $X:=Z_{\R}(g^*+\epsilon P) \subset \mathbb{RP}^4$ is smooth for all $\epsilon \in [0,1)$. 
\end{proposition}

\begin{proof}
    The condition $\|{\tilde{P}}_{C^0} < k$ ensures that $\tilde{g}+\epsilon \tilde{P}$ has no zeros outside of $U$, i.e. in particular the variety $Z(\tilde{g}+\epsilon \tilde{P})$ has no singular points outside of $U$.
    The condition $\|{\nabla \tilde{P}}_{C^0} < m$ ensures that $\nabla (\tilde{g}+\epsilon \tilde{P}) \neq 0$ on $U$, in particular the variety $Z(\tilde{g}+\epsilon \tilde{P})$ has no singular points in $U$.
\end{proof}

A choice for $U$ that is useful for computations, because numerical minimisation on it works reliably, is $U:=\{x \in S^4: |g^*(x)| < \rho \}$ for some fixed $\rho>0$.
For $P:=x_0^5+x_1^5+x_2^5+x_3^5+x_4^5$, one checks that
\begin{align}
    \texttt{Quintic}
    :=
    Z(g^*+ \epsilon P)
    \subset \mathbb{CP}^4
    \text{ for }
    \epsilon = \frac{27}{1000},
\end{align}
i.e. $g^*+ \epsilon P=(x_0^2+\dots+x_4^2)\left( x_0(x_1^2+x_2^2+x_3^2-x_4^2)-(x_1^3+x_2^3+x_3^3-
\frac{1}{2}x_4^3-\frac{1}{4}x_0^3 \right)+\frac{27}{1000}(x_0^5+x_1^5+x_2^5+x_3^5+x_4^5)$,
is smooth.

Using numerical minimisation together with \cref{proposition:conrete-small-deformation-bound} suggests that also $Z(g^*+ \epsilon P)$ is smooth for every $\epsilon \leq \frac{27}{1000}$, however we remark that due to the use of numerical minimisation we have no rigorous proof of that fact.
Because of this, we cannot guarantee that the real locus is diffeomorphic to $Z(f) \subset \mathbb{RP}^4$.
Using the persistent homology tool from \cite{Gafvert2022,DiRocco2022} suggests that $b_1(Z(g^*+\epsilon P)=1$, though also here we did no provably correct calculation of this Betti number (which would be available using the same tool but is computationally very expensive).

\subsubsection{Intersections of a quadric and a quartic \texttt{CICY1} and \texttt{CICY2} with $b_1(L)=1$}
\label{subsection:intersections-of-quadrics-and-quartics}

Let $Q_{\text{aff}}(x_1,x_2,x_3)$ be a quartic defining a smooth variety in $\mathbb{C}^3$ whose real locus in $\mathbb{R}^3$ is compact and non-empty, e.g. the deformed sphere $Q_{\text{aff}}(x_1,x_2,x_3)=x_1^4+x_2^4+x_3^4-1$.
Let $C_{\text{aff}}(x_4,x_5)=x_4^2+x_5^2-1$ be the equation of a circle in $\mathbb{R}^2$.

Then
\[
X^*_{\text{aff}}:=Z_{\R}(Q_{\text{aff}}(x_1,x_2,x_3), C_{\text{aff}}(x_4,x_5)) \subset \mathbb{R}^5
\]
is diffeomorphic to $S^1 \times Z(Q_{\text{aff}})$, and therefore has $b_1(X^*_{\text{aff}})>0$.
Denote the projectivisation in $\mathbb{CP}^5=[z_0:\dots:z_5]$ of $Q_{\text{aff}}$ and $C_{\text{aff}}$ by $Q$ and $C$ respectively.
Then $X^*:=Z(Q, C) \subset \mathbb{CP}^5$ has a one-dimensional singularity at $Z(x_1^4+x_2^4+x_3^4, \, x_0, \, x_4, \, x_5) \subset \mathbb{CP}^5$ and two point singularities at
\[Z(x_4^2+x_5^2, \, x_0, \, x_1, \, x_2, \, x_3) = \{ [0:0:0:0:1:i], [0:0:0:0:1:-i] \} \subset \mathbb{CP}^5.\]

The singularities are away from the real locus, so $X^* \cap \RP^5=Z_{\R}(Q,C) \subset \mathbb{RP}^5$ is smooth as a real variety.
Therefore, we can study small perturbations of $X^*$ as in the previous section.
Let $P_k(x_0,\dots,x_5)$ be a generic homogeneous polynomial of degree $k$.
Then, for small $\epsilon > 0$, we have that $X_\epsilon:=Z(Q+\epsilon P_4, C+\epsilon P_2) \subset \mathbb{CP}^5$ is smooth (as a complex variety), and the real loci $X_\epsilon \cap \RP^5$ and $X^* \cap \RP^5$ are diffeomorphic for small $\epsilon$, in particular $b_1(X \cap \RP^5) > 0$.
We can be more explicit:

\begin{proposition}
    \label{proposition:smooth-deformation-of-cicy}
    Let $P_2=x_1^2+x_2^2$ and $P_4=x_4^4$.
    Let $\epsilon_0 \approx 0.590$ be the unique real root of $1-x-2x^3=0$.
    Then, for $0<\epsilon < \epsilon_0$, the variety $X_\epsilon$ is smooth.
\end{proposition}

\begin{proof}
    Fix $0<\epsilon < \epsilon_0$ and denote $F=(Q+\epsilon P_4, C+\epsilon P_2) : \mathbb{C}^6 \rightarrow \mathbb{C}^2$.
    Then
    \[
        DF = 
        \begin{pmatrix}
            -4x_0^3 & 4x_1^3 & 4x_2^3 & 4x_3^3 & 4\epsilon x_4^3 & 0
            \\
            -2x_0 & 2\epsilon x_1 & 2\epsilon x_2 & 0 & 2x_4 & 2x_5
        \end{pmatrix}
    \]
    and
    \[
        \sing(X_\epsilon)
        =
        \{
        x \in \mathbb{CP}^5 :
        F(x)=0 \text{ and }
        \rank(DF)(x)<2
        \},
    \]
    and we will show that $\sing(X_\epsilon) = \emptyset$.
    Assume there exists $x \in \sing(X_\epsilon)$.
    Because $\rank(DF)(x)<2$ we either have that one of the rows of $\rank(DF)(x)$ is zero or that the two rows are a non-zero multiple of each other.
    The first case gives a contradiction, because there is no point on $X_\epsilon$ which has five coordinates set to zero.
    The second case implies $x_3=x_5=0$.
    
    We distinguish the following cases and will obtain a contradiction in each of them:

    \textbf{Case I: $x_0 \neq 0$.}
    Without loss of generality $x_0=1$.
    Then
    \begin{align}
    \label{equation:smoothness-proof-coordinate-expressions}
    \begin{split}
        \epsilon x_1 &= x_1^3
        \Rightarrow
        x_1 \in 
        \left\{
        0, \pm \sqrt{\epsilon} 
        \right\},
        \\
        \epsilon x_2 &= x_2^3
        \Rightarrow
        x_2 \in 
        \left\{
        0, \pm \sqrt{\epsilon} 
        \right\},
        \\
        x_4 &= \epsilon x_4^3
        \Rightarrow
        x_4 \in 
        \left\{
        0, \pm \sqrt{\frac{1}{\epsilon}} 
        \right\}.
    \end{split}
    \end{align}
    
    \textbf{Case I.A: $x_4=0$.}
    We have the following due to $(C+\epsilon P_2)(x)=0$:
    \begin{align*}
        0
        &=
        |-1+\epsilon x_1^2+\epsilon x_2^2+x_4^2+x_5^2|
        \geq
        1-\epsilon |x_1|^2-\epsilon |x_2|^2
        \geq
        1-2\epsilon^2
    \end{align*}
    which is $>0$ because $\epsilon<\epsilon_0<\frac{\sqrt{2}}{2}$.
    Here we used the triangle inequality for the first inequality and \cref{equation:smoothness-proof-coordinate-expressions} for the second inequality.
    This is a contradiction.
    
    \textbf{Case I.B: $x_4 \neq 0$.}
    In this case,
    \begin{align*}
        0
        &=
        |-1+\epsilon x_1^2+\epsilon x_2^2+x_4^2+x_5^2|
        \geq
        |x_4|^2-1-2 \epsilon^2
        =
        \frac{1}{\epsilon}-1-2 \epsilon^2
    \end{align*}
    which is $>0$ by the assumption that $\epsilon < \epsilon_0$ and the definition of $\epsilon_0$.
    This is again a contradiction.

    \textbf{Case II: $x_0 = 0$. Case II.A: $x_1 \neq 0$.}
    Without loss of generality set $x_1=1$.
    Similar to before, we get 
    \begin{align}
    \label{equation:smoothness-proof-coordinate-expressions-case-II}
    \begin{split}
        x_2 &= x_2^3
        \Rightarrow
        x_2 \in 
        \left\{
        0, \pm 1 
        \right\},
        \\
        x_4 &= \epsilon^2 x_4^3
        \Rightarrow
        x_4 \in 
        \left\{
        0, \pm \frac{1}{\epsilon} 
        \right\}.
    \end{split}
    \end{align}
    In particular, $x_1^4+x_2^4+x_3^4+\epsilon x_4^4$ is the sum of non-negative real numbers, so $Q(x)=0$ implies $x_1=0$ which is a contradiction.

    \textbf{Case II.B: $x_1 = 0 $:}
    The equations $x_3=x_5=0$ imply $x_2 \neq 0$ and $x_4 \neq 0$.
    Rescale to $x_2=1$, then as before $x_4 \in \{0, \pm 1/\epsilon \}$, and we obtain a contradiction as in case II.A.
\end{proof}

Note that the singular variety $X^*$ has a $\text{U}(1)$-symmetry, namely a rotation in the $(x_4,x_5)$-plane.
If $X^*$ was smooth, it would admit a Calabi-Yau metric, the infinitesimal $\text{U}(1)$-action would define a Killing field, and its dual would be a parallel $1$-form due to $\Ric=0$ and the Bochner trick (see \cite{Bochner1946} or \cite[Theorem 3.1]{Wu2017}).
Any smoothing of $X^*$ must lose this $\text{U}(1)$-symmetry, as a compact Calabi-Yau manifold with holonomy $\text{SU}(m)$ cannot have a continuous symmetry.
However, one may speculate that for small values of $\epsilon$, the infinitesimal $\text{U}(1)$-action projected onto $X_\epsilon$ gives rise to a $1$-form that is close to being parallel, and that there is a nearby harmonic nowhere vanishing $1$-form.
In \cref{section:experimental-results} we present numerical evidence suggesting that this is indeed the case.

We study the following two varieties:

\begin{align}
    \begin{split}
    \text{\texttt{CICY1}}
    =
    Z
    \left(
    -x_0^4+x_1^4+x_2^4+x_3^4+\frac{1}{4}x_4^4,
    \right.
    &
    \left.
    -x_0^2+\frac{1}{4}x_1^2+\frac{1}{4}x_2^2+x_4^2+x_5^2
    \right),
    \\
    \text{\texttt{CICY2}}
    =
    Z
    \left(
    -x_0^4+x_1^4+x_2^4+x_3^4+\frac{1}{10}x_4^4+\frac{1}{10}x_5^4,
    \right.
    &
    \left.
    -x_0^2+\frac{1}{100}x_1^2+\frac{1}{100}x_2^2+\frac{1}{100}x_3^2+x_4^2+x_5^2
    \right).
    \end{split}
\end{align}

By \cref{proposition:smooth-deformation-of-cicy}, \texttt{CICY1} is smooth and its real locus is diffeomorphic to $S^1 \times S^2$.
The defining equations are far away from the singular limit $X^*$, so we can say nothing about harmonic $1$-forms on its real locus.

One can adapt the proof of \cref{proposition:smooth-deformation-of-cicy} to show that also \texttt{CICY2} is smooth and its real locus is diffeomorphic to $S^1 \times S^2$.
The defining equations are close to the singular limit $X^*$, and because of the discussion of the $\text{U}(1)$-symmetry above we speculate that its real locus admits a nowhere vanishing harmonic $1$-form.

We will pick this discussion up again in \cref{section:experimental-results}.

\section{Numerics}

In this section we describe our numerical experiments.
Our implementations of these algorithms can be found at \url{https://github.com/yidiq7/MLHarmonic-1-form}.

\subsection{ Sampling algorithms }
\label{subsection:sampling-algorithms}

The algorithms in the following subsections require doing numerical integrals over a K\"ahler manifold $M$
or its real locus $L$.
For $f:M \rightarrow \R$ and a measure $\d\mu_S$ on $M$ consider the integral 
\[
\bar f = \int_M \d\mu_S \, f .
\]
We will estimate it by averaging over a subset $S_N\subset M$ of $N$ independent samples from a probability measure $\d\mu_S$.
Define the random variable
\[
\CI_{S,N}[ f ] \equiv \frac{1}{N}\sum_{x\in S_N} f(x) ,
\]
and let $\E[ \cdot ]$ be the expectation under sampling from $\d\mu_S$, then
\[
\E\bigg[\CI_{S,N}[f]\bigg] \equiv \bar f = \int_M \d\mu_S \, f
\]
and the error introduced by the sampling is of the order of the square root of
\begin{align*}
\Var\bigg[\CI_{S,N}[f]\bigg] 
&\equiv \E\bigg[\CI_N[f]^2\bigg] - \E\bigg[\CI_N[f]\bigg]^2 
= \frac{1}{N} \int \d\mu_S (f - \bar f)^2 .
\end{align*}
More generally, we can take the sampling distribution $\d\mu_S$ different from the Kähler volume form $\d\mu$.
To minimize the sampling variance of $\int \d\mu\, f$, 
we want to choose $\d\mu_S$ to minimize the variance of $f\,\d\mu/\d\mu_S$.
Ideally one would take $\d\mu_S = f \d\mu$ in which case one can take $N=1$.
This requires knowing $f$ exactly, which we don't in our problems.
However, the function $f$ is bounded and so the main point is to control the ratio of the volume forms $\d\mu/\d\mu_S$.
To do this it is good enough to take $\d\mu_S$ to be the Fubini-Study metric induced from an embedding into projective space.

This can be accomplished as follows.
In \cite{shiffman_distribution_1999}, the distribution of the simultaneous zeroes of a $d$-tuple of random sections on $M^d$ was found.
The idea is to use the Poincar\'e-Lelong formula for the $(1,1)$-form representing the divisor (zero set),
\[
Z_s = \frac{i}{\pi}\partial\bar\partial \log ||s|| ,
\]
and integrate over the choice of $s$.
The random sections are drawn from a normal distribution characterized by a hermitian form $G$.  If these are a complete set of sections of an ample
bundle $L$, there is an embedding $\iota:M^d\hookrightarrow \IP H^0(L)$ and $G$ defines a Fubini-Study metric on the ambient space.  
Let its K\"ahler 
form be $\omega_G$, then
\be\label{eq:shifzel}
\mu_S = (\iota^*\omega_{FS})^d .
\ee
This can be used to sample from $M^d$ which is a complete intersection in $\IC\IP^N$.  The direct approach is to choose $d$ random
sections on $\IC\IP^N$ and intersect their zeroes with those of the defining sections.  
One can instead dualise and choose $N+1-d$ random points on $\IC\IP^N$, whose span is a $\IC\IP^{N-d}$, and take the zeroes
of the defining sections on this span.

The method can be used for any ambient space which can be embedded by an explicit complete basis of sections.
It was generalised to hypersurfaces in toric varieties in \cite{Braun2008}.
The challenge there is to work out the map between the toric coordinates and the sections;
methods for this are discussed in \cite{Larfors2022}.

Next we want to sample points on the real locus $L\subset M$.
Note that, up to conjugation, the complex conjugation $z \mapsto \bar{z}$ is the only anti-holomorphic involution on $\mathbb{CP}^N$ for even $N$, see \cite[Lemma 11]{Morozov2022}, so there is little loss of generality in taking the real locus rather than the fixed set of an arbitrary anti-holomorphic involution.
A sampling algorithm for real algebraic submanifolds using the technique of intersection with a linear space was
developed in \cite{Breiding2020} (see also their references \cite{Boykov2003,Legland2007,Lehmann2012,Li2003}).
The result most relevant for us is \cite[Theorem 6.1]{Breiding2020}, which we now adapt to our conventions.

Let $L$ be an algebraic variety in $\IR\IP^N$ and $d=\dim L$.  
Let $A\in \IR^{d\times(N+1)}$ be a random matrix with
i.i.d. entries in $\CN(0,1)$.  Define the linear space $\CL_A \equiv \{ x\in\IR\IP^N | Ax=0 \}$.  Then
\begin{align*}
\int_L \d\mu_{FS}\, f &= \mbox{Vol}(\IR\IP^N)\, \E[{\bar f}(A)] \\
{\bar f}(A) &\equiv \sum_{x\in L\cap \CL_A} f(x) .
\end{align*}
One can then approximate this integral by sampling $A$ and averaging over the subsets $x\in L\cap \CL_A$.
In words, these are the intersections of linear subspaces of $\IR\IP^N$ with $L$, so this prescription also expresses the measure in terms of zeroes of random sections (here defined
by the rows of $A$).
One can also sample the same measure by taking $\CL_A$ to be the linear span in $\IR^{N+1}$ of $N+1-d$ points sampled from $\d\mu_{FS}$.
This follows because there is a unique $SO(N+1)$-invariant probability measure on the Grassmannian manifold of subspaces $\CL$.

To integrate over a manifold embedded in a product of $\IR\IP^N$'s, one could iterate this procedure.
The result would be an expectation of a sum of points $x\in L \cap L_A \cap L_B \ldots \cap L_X$, where each subspace $L_A$, $L_B$, $\dots$ is sampled from one of the projective spaces.

\medskip

\subsection{Approximating the Calabi-Yau metric using bihomogeneous neural networks}
\label{section:bihomogeneous-nn}

As announced in \cref{subsection:construct-nowhere-vanishing-from-approximate}, we will compute approximate Calabi-Yau metrics on several Calabi-Yau varieties.
This problem has a rich history:
it was first studied in \cite{Headrick2005} and subsequently \cite{Douglas2008,Braun2008}.
In \cite{Donaldson2009} a new approach based on embedding the Calabi-Yau variety into a high-dimensional projective space was pioneered, and this idea still underlies many of the modern approaches.
In the subsequent work \cite{Headrick2013} the numerics were improved, and this approach currently provides the best approximate Calabi-Yau metrics on varieties with many discrete symmetries.
For varieties without symmetry, the machine learning approaches \cite{Ashmore2020,Anderson2021,Douglas2022,Douglas2022b} as well as \cite{jejjala_neural_2022, larfors2021learning, Berglund2022,gerdes_cyjax_2023} emerged, and appear to provide the best approximate metrics on such varieties.
We also note \cite{Doran2008} as well as \cite{Bunch2008,Seyyedali2009,Anderson2010} on very closely related numerical problems.

First, to do a numerically verified version of the analysis in \cref{subsection:construct-nowhere-vanishing-from-approximate} it is important that the approximate metric be differentiable.
(We do not do this analysis in this paper, but hope to do so in the future.)
Second, our varieties have few discrete symmetries.
These two features restrict the methods we can choose from.
We decided to use the \emph{bihomogeneous neural networks} from \cite{Douglas2022}.
An implementation for hypersurfaces can be found at \cite{QiGithub2024}, and we expanded this code to work for complete intersection Calabi-Yaus in a single projective space.
A detailed description of the approach can be found in \cite[Section 2.4]{Douglas2022}, and we briefly review it here.

Consider a Calabi-Yau manifold $Y \subset \mathbb{CP}^n$ defined as a complete intersection of polynomials, such as \texttt{Fermat} or \texttt{CICY1}, and denote coordinates on $\mathbb{CP}^n$ by $z_0,\dots,z_n$.
We will now define a feed-forward neural network
\begin{align}\label{eq:bihomo-net}
    F(z \bar{z})
    =
    W^{(d)} \circ \theta \circ W^{(d-1)} \circ \dots \circ \theta \circ W^{(1)} \circ \theta \circ W^{(0)}(z \bar{z}),
\end{align}
where the $W^{(i)}$'s are real general linear transformations, and $\theta$ is the activation function defined in every component as $z \mapsto z^2$.
Here, $z \bar{z}$ denotes the real and imaginary parts of all bihomogeneous combinations $z_i \bar{z}_j$, which in complex dimension $n$ are $n^2$ inputs.
Note that in this way, $F$ is \emph{not} a well-defined function on $\mathbb{CP}^n$, but the Kähler form $\frac{i}{2} \partial \bar{\partial} \log(F)$ is well-defined.
That is in analogy to the Kähler form of the ordinary Fubini-Study metric, which is $\omega_{\text{FS}}=\frac{i}{2} \partial \bar{\partial} \log \left( \sum_i z_i \bar{z}_i \right)$.

Fixing $d \in \mathbb{N}$ and the dimensions of each $W^{(i)}$ for $i=1,\dots,d$, the ideal goal would be to find entries of the matrices $W^{(i)}$ so that the Kähler form $\omega = \frac{i}{2} \partial \bar{\partial} \log \left( F(z \bar{z}) \right)$ solves \cref{equation:cy-kaehler-and-holo-volume-form-comparison}, i.e. its Riemannian volume form is $\Omega \wedge \bar{\Omega}$.
It is in general impossible to find matrix entries to solve the equation exactly, so instead one aims to find matrix entries to minimise the \emph{Monge–Ampère loss functional}
\begin{align}\label{eq:MSE-loss}
\int_Y
\left( \frac{\omega^3}{\Omega \wedge \bar{\Omega}}-1 \right)^2
\operatorname{vol},
\end{align}
where $\operatorname{vol}$ is any fixed volume form on $Y$.
Because it is practically impossible to exactly evaluate this integral, we use Monte-Carlo integration over the point samples explained in \cref{subsection:sampling-algorithms}.

This is a minimisation problem on a very high-dimensional space, typically of $100,000$ parameters or more.
As for exactly computing an integral over $Y$, it is also practically impossible to find a global minimiser of the Monge–Ampère loss.
Instead, as is the standard thing to do in this situation, we use gradient descent to find a local minimiser.
We achieve this by using the implementation of gradient descent from the machine learning package \texttt{tensorflow}.

\subsection{Numerical forms using a polynomial basis}
\label{section:numerical-forms}

Since our surface $L$ is an algebraic hypersurface in real projective space, we can use algebraic techniques to
represent and solve for the harmonic one-form.  This has two advantages.  
First, the harmonic form is smooth so a polynomial basis should be highly convergent.
Second, it is similar to our existing code for Calabi-Yau metrics and easy to implement.

Thus, we want to find a harmonic form on $L=Z_{\mathbb{R}}(f) \subset \IR\IP^4$, where the metric $g$
on $L$ is the restriction of a given K\"ahler metric on the complex hypersurface. (In this section we will work out the details for quintic hypersurfaces. The calculations for CICYs can be done in the same manner on $\IR\IP^5$.)

A harmonic form $\omega$ minimizes a quadratic energy on $L$,
\be \label{eq:quadE}
E = \int_L |\d \omega|^2 + |\d *_L\omega|^2 ,
\ee
subject to the constraint
\be \label{eq:norm-constraint}
1 = \int_L |\omega|^2 .
\ee
Here $|\cdot|$ is the pointwise norm on $L$ using the metric $g$  which can also be written $|\alpha|^2=(\alpha\wedge *_L\alpha)/ \operatorname{vol}$.
Given a parameterized ansatz $\omega_t \in \Omega^1(L)$ we can get an upper bound on this by minimizing $E(t)$.

We have a projection map onto $\langle \grad f \rangle^\perp$, say $P_L: T \mathbb{RP}^4 \rightarrow T \mathbb{RP}^4$, which induces projections on all differential forms that we denote by the same symbol by abuse of notation.
Given a parameterized ansatz of forms $\omega_t \in \Omega^1(\IR\IP^4)$, the equations \cref{eq:quadE,eq:norm-constraint} for the form $P_L(\omega_t)|_L$ become:
\begin{align*}
E_1 &= \int_L |P_L\, \d \,P_L\omega_t|^2 + |P_L\, \d \,*_L P_L\omega_t|^2, \\
1 &= \int_L |P_L\omega_t|^2.
\end{align*}

The integrals can be estimated by sampling on $L$ and the constraint can be imposed by normalizing $\omega$ in each gradient descent step.
Thus the remaining work is to choose a basis and implement $P_L$, $*_L$ and the norm $|\cdot |$.

We first define a set of globally defined one-forms $e^{[ij]}$, $1\le i<j\le 5$, on $\IR\IP^4$ as the basis.
We define $\IR\IP^N$ as the quotient of $\IR^{N+1}-\{0\}$ by the action $x^i\rightarrow \lambda x^i$,
generated by the vector field $\vec e = \sum_{i=1}^{N+1} x^i \partial/\partial x^i$.
A geometric object $\alpha$ on $\IR^{N+1}$ will descend to the quotient if
$\mbox{Lie}_{\vec e} \alpha = 0$.
Generally this is solved by taking $\alpha$ homogeneous.
In particular we can define functions $u^i=x^i/|x|$ on $S^N$ and
$u^{ij} = x^i x^j/|x|^2$ on $\IR\IP^N$.  These satisfy $\sum_i u^{ii}=1$.
They satisfy many more relations of course such as $u_{ij}u_{kl}=u_{ik}u_{jl}$.

Next, $T_x\IR\IP^N$ is the quotient of $T_x \IR^{N+1}$ by the normal bundle $N_x$ spanned by $\vec e_x$.
Its dual $T^*_x \IR\IP^N$ is the subspace of $T^*_x \IR^{N+1}$ annihilated by $\vec e_x$,
so well defined forms on $\IR\IP^N$ must be annihilated by $\vec e_x$.
Given a one-form $\lambda\in T^* \IR^{N+1}$ this can be accomplished by taking $\lambda - (\vec e \lambda) E$
where $E \equiv \sum_i x^i \d x^i/|x|^2$.

Using this idea we can define an overcomplete basis of one-forms 
\[
e^{ij} = \frac{ x^i }{ |x|^2 } \d x^j - \frac{x^i x^j}{ |x|^2} E .
\]
Denote the symmetric and antisymmetric parts as
\begin{align}
\notag
e^{(ij)} &= \frac{1}{2}\left(e_{ij}+ e_{ji}\right) \\
e^{[ij]} &= \frac{1}{2}\left(e_{ij} -e_{ji}\right) = \frac{ x^i \d x^j - x^j \d x^i }{ 2|x|^2 } .  \label{eq:e-basis}
\end{align}
The $e^{[ij]}$ are an overcomplete basis, {\it i.e.} they span $T^*_x \IR\IP^N$ for all $x \in \IR\IP^N$.
To see this, note that at $x_{\text{base}}:=[1:0:\ldots:0]$ the $e^{[1i]}$ with $i>1$ form a basis.
The $SO(N+1)$-action then maps the overcomplete basis from $x_{\text{base}}$ to any other point in $\IR\IP^N$.
Furthermore 
$ \d \, (1/|x|^2) = -2 E/|x|^2 $
and
\begin{align*}
\d u^{ij} &= (x_i \d x_j + x_j \d x_j)/|x|^2 - x_i x_j E / |x|^2 = 2 e^{(ij)} \\
\d e^{(ij)} &= 0 \\ 
\d e^{[ij]} &= \frac{\d x^i \wedge \d x^j }{ |x|^2 } - 2 E \wedge e^{[ij]}
= \sum_k e^{ki} \wedge e^{kj} .
\end{align*}
This last equation is reminiscent of the Cartan structural equations and if we choose a frame could be used to derive them.
Now we can express $\omega$ as
\be\label{eq:defomega}
\omega = \sum_{1\le i<j\le 5} \Omega_{ij} e^{[ij]},
\ee

where $\Omega_{ij}$ are polynomials represented by the outputs of a homogeneous neural network
\be
    F(x)
    =
    W^{(d)} \circ \theta \circ W^{(d-1)} \circ \dots \circ \theta \circ W^{(1)} \circ \theta \circ W^{(0)}(x / |x|) ; \quad x \in \IR\IP^4
\ee
with the activation function $\theta$ defined as $x \mapsto x^2$.
Here, $|x|$ denotes the $l^2$-norm on $\IR^5$.
This is the same network architecture as in \cref{eq:bihomo-net}, only now all inputs are \emph{real}, so we need not consider bihomogeneous combinations of input coordinates.

After plugging in Eq. \ref{eq:e-basis}, the form can be written in the basis of $\d x^k$ as follows:
\be
\omega = \sum_{1\le i<j\le 5} \Omega_{ij} \frac{ x^i \d x^j - x^j \d x^i }{ 2|x|^2 } = \omega_k \d x^k.
\ee
The simplest way to proceed is to restrict the forms to local coordinates on $L$ and use these to do all the further computations.
In this case we can use Eq. \ref{eq:quadE} directly and do not need further use of $P_L$.
We define patches $U_i\subset \IC\IP^4$ with $z^i\ne 0$ and $U_{i,j} \subset M = Z(f) \subset \mathbb{CP}^4$ in which we solve the equation $f=0$ for $z^j$.
(In the CICY case of $Z(f_1,f_2) \subset \IC\IP^5$ we solve the system of equations $f_1 = f_2 = 0$ for for $z^j$.)
The local coordinates on $U_{i,j}$ are defined by setting $z^i=1$, solving for $z^j$ and taking the three remaining coordinates.
These patches can also be used for $L$ with coordinates $x^i=\mbox{Re}\, z^i$.  The form $\d x^i=0$ in $U_i$, and 
to get $\d x^j$ in $U_{i,j}$ one solves $\d f=0$, so $\d x^j = -\sum_k (\partial_k f/\partial_j f) \d x^k$.

Now the first half of the functional $E$ can be computed as
\be \label{eq:d-omega-int}
\int |\d\omega|^2 = \int (\det g)^{1/2} \left(\partial_{[l} \omega_{k]} g^{km} g^{ln} \partial_{[n} \omega_{m]} \right),
\ee
where $g$ is the numerical metric computed in Section \ref{section:bihomogeneous-nn}
and the partial derivatives can be computed using auto-diff, then be made antisymmetric manually.

It remains to implement $*_L$ acting on one-forms and two-forms.
The coordinate expression is
\[
*_L \d x^k = \sqrt{g} \epsilon^{klm} g_{lp} g_{mq} \d x^p \wedge \d x^q.
\]
If we are only using this to write the norm and $* \d *\omega$, a bit simpler is to write this as the divergence 
\[
\operatorname{div} \left( \sum_k\omega_k \d x^k \right) = g^{kl} \nabla_l \omega_k = (\det g)^{-1/2} \partial_l \left((\det g)^{1/2} g^{kl} \omega_k \right)
\]
so 
\[
\int |\d *\omega|^2 = \int (\det g)^{-1/2} \left(\partial_l \left((\det g)^{1/2} g^{kl} \omega_k \right)\right)^2 .
\]
In practice, since $g$ and $\omega$ are represented by neural networks, it would be computationally too expensive to trace both of them simultaneously in auto-diff during the training process. 
Therefore we choose to rewrite this integral as
\be \label{eq:d-star-omega-int}
\int |\d *\omega|^2 = \int \left(\det g \right)^{-1/2} \left(  \partial_l \left(\det g \right)^{1/2}  g^{kl} \omega_k + \left(\det g\right)^{1/2}\left(g^{kl} \partial_l \omega_k  + \omega_k \partial_l g^{kl}\right) \right)^2.
\ee

Combining \ref{eq:d-omega-int} and \ref{eq:d-star-omega-int}, we have the full expression for the functional $E$. We will use it as the loss function to train the network. The approximate harmonic 1-form will be obtained once $E$ is minimized.

\section{Experimental Results}
\label{section:experimental-results}

We conduct experiments on $4$ Calabi-Yau manifolds:
\begin{itemize}
    \item 
    The Fermat quintic \texttt{Fermat} from \cref{section:fermat},

    \item 
    The quintic \texttt{Quintic} from \cref{section:quintic},

    \item 
    The intersections of a quadric and quartic \texttt{CICY1} and \texttt{CICY2} from \cref{subsection:intersections-of-quadrics-and-quartics}.
\end{itemize}
Our experiments consist of two steps:
\begin{enumerate}
    \item 
    Compute an approximation $g$ of the Calabi-Yau metric $g_{\text{CY}}$ using the technique explained in \cref{section:bihomogeneous-nn}.
    The results of this process are presented in \cref{subsection:results-approximate-cy-metrics}.

    \item 
    Approximately compute a $1$-form $\lambda \in \Omega^1(L)$ that is harmonic with respect to $g \mid_L$ using the technique from \cref{section:numerical-forms}.
    The results of this are presented in \cref{section:results-1-forms}.
\end{enumerate}

\subsection{Approximate Calabi-Yau metrics}
\label{subsection:results-approximate-cy-metrics}

\paragraph{Smoothness and approximation error.}

The error of solving the complex Monge–Ampère equation has been called Monge–Ampère loss in the literature.
It has been observed in \cite[Figures 3 and 4]{Douglas2022}, \cite[Figure 12]{Anderson2021} and \cite[Figure 10]{Headrick2013} that the quality of the approximate Calabi-Yau metric is related to how close the Calabi-Yau manifold is to being singular.
Compare also \cite[Section 7.3]{Li2023} discussing this.
We show the Monge–Ampère loss for the four varieties \texttt{Fermat}, \texttt{CICY1}, \texttt{CICY2}, and \texttt{Quintic} in \cref{tab:sigma_loss} confirming this observation.

\begin{table}[h]
    \centering
    \begin{tabular}{|c|c|}
        \hline
        Manifold & Monge–Ampère loss\\
        \hline
        \texttt{Fermat} & $8.16 \times 10^{-4}$\\
        \texttt{Quintic} & $7.13 \times 10^{-2}$\\
        \texttt{CICY1} & $2.56 \times 10^{-3}$\\
        \texttt{CICY2} & $5.68 \times 10^{-3}$\\
        \hline
    \end{tabular}
    \caption{The Monge–Ampère loss measures how far an approximate Calabi-Yau metric is from solving the complex Monge–Ampère equation. We obtain the best approximation on the very smooth variety \texttt{Fermat}. Comparing \texttt{CICY1} and \texttt{CICY2}, we observe that the variety that is closer to the singular limit has larger Monge–Ampère loss.}
    \label{tab:sigma_loss}
\end{table}

\paragraph{Approximation error in different regions.}

\begin{figure}[htbp]
    \centering
    \begin{minipage}{0.24 \textwidth}
        \centering
        \includegraphics[width=4cm]{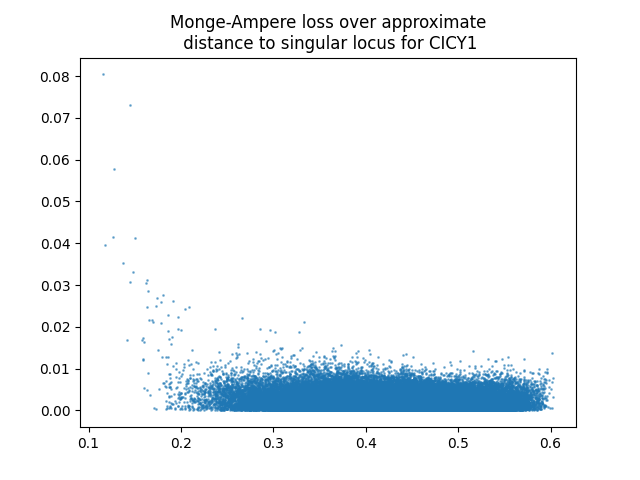}
    \end{minipage}
    \begin{minipage}{0.24 \textwidth}
        \centering
        \includegraphics[width=4cm]{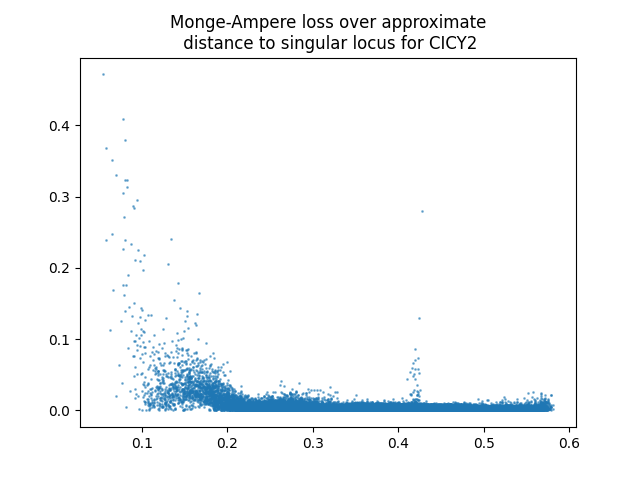}
    \end{minipage}
    \begin{minipage}{0.24 \textwidth}
        \centering
        \includegraphics[width=4cm]{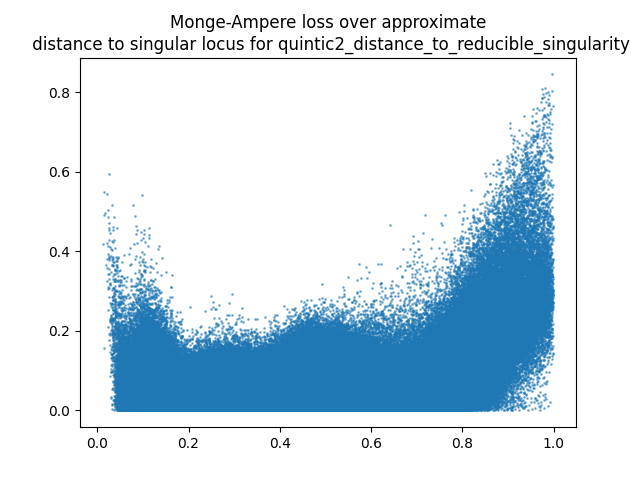}
    \end{minipage}
    \begin{minipage}{0.24 \textwidth}
        \centering
        \includegraphics[width=4cm]{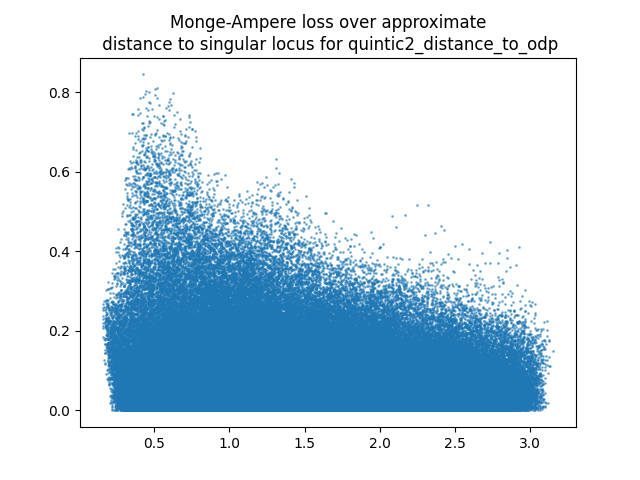}
    \end{minipage}
    \caption{Plots showing the Monge–Ampère error $|1-\omega^m/\Omega \wedge \overline{\Omega}|$ on the $y$-axis over the approximate distance to a singular locus. From left to right: \texttt{CICY1}, \texttt{CICY2}, then two plots for \texttt{Quintic}. The approximate distance to a singular locus is computed as follows:
    for \texttt{CICY1} and \texttt{CICY2} as $\min \{d_2(x)/|x|^2, \, d_4(x)/|x|^4\}$, where $d_2(x)=\max \{ |x_4^2+x_5^2|, \, |x_0|^2, \, |x_1|^2, \, |x_2|^2, \, |x_3|^2 \}$ and $d_4(x)=\max \{ |x_1^4+x_2^4+x_3^4|, \, |x_0|^4, \, |x_4|^4, \, |x_5|^4 \}$.
    \texttt{Quintic} is a smoothing of a variety with two singularities, namely (a) the intersection of the divisors $Z(f)$ and $Z(g^*)$ from \cref{section:quintic}, and (b) an ordinary double point at $a=[1:0:0:0:0]$.
    The $x$-axis of the third image is an approximation of the distance to (a), namely $\max \{|x_0^2+\dots+x_4^2|/|x|^2, \, f^*(x)/|x|^3\}$.
    The $x$-axis of the fourth image is an approximation of the distance to (b), namely $|x-(|x|,0,0,0,0)|/|x|$.}
    \label{figure:metric-approximation-quality}
\end{figure}

Apart from this, it is folklore that for a Calabi-Yau manifold that is a smoothing of a singular Calabi-Yau manifold, the approximate metric is only a poor fit near the smoothed singularity, and is a good fit far away from the smoothed singularity.
This is confirmed by the data presented in \cref{figure:metric-approximation-quality}.

\paragraph{Formation of long necks on Calabi-Yau manifolds that are smoothings of singularities.}

The variety \texttt{Quintic} is defined as the perturbation $Q$ of a singular variety $Q^*$, see \cref{section:quintic}.
It is known from \cite[Theorem 1.2]{Sun2019} that as the perturbation $Q$ approaches $Q^*$, a long neck must be forming near the singular locus of $Q^*$.

Likewise, the variety \texttt{CICY2} is defined as the perturbation of s singular variety $X^*$, see \cref{subsection:intersections-of-quadrics-and-quartics}.
The work \cite{Sun2019} applies to hypersurfaces and not codimension two subvarieties such as \texttt{CICY2}.
Nevertheless, we observe experimentally that a long neck is forming as the smoothing of the singular Calabi-Yau approaches the singular limit $X^*$.

\begin{figure}[htbp]
    \centering
    \begin{minipage}{0.3 \textwidth}
        \centering
        \includegraphics[width=5cm]{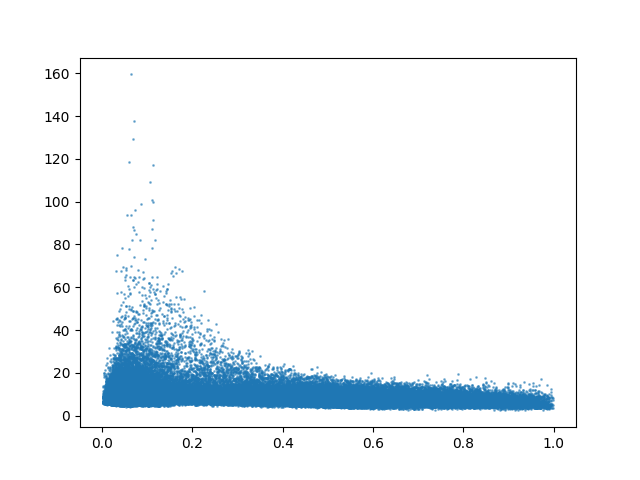}
    \end{minipage}
    \begin{minipage}{0.3 \textwidth}
        \centering
        \includegraphics[width=5cm]{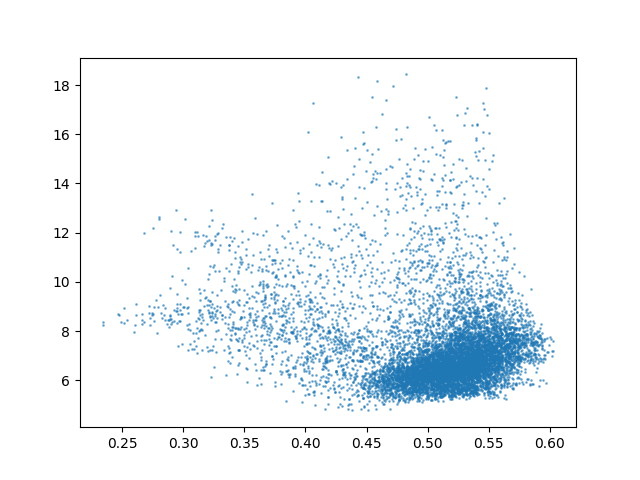}
    \end{minipage}
    \begin{minipage}{0.3 \textwidth}
        \centering
        \includegraphics[width=5cm]{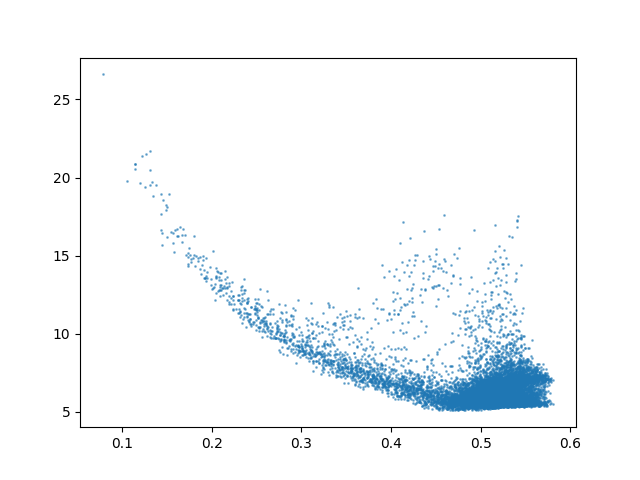}
    \end{minipage}
    \caption{Plots showing the formation of long necks. From left to right: \texttt{Quintic}, \texttt{CICY1}, \texttt{CICY2}. The three varieties are perturbations of singular varieties. The horizontal axis measures an approximation of the distance to the singular locus of the initial singular variety, namely $\max \{|x_0^2+\dots+x_4^2|/|x|^2, \, f^*(x)/|x|^3\}$ for \texttt{Quintic} and $\min \{d_2(x)/|x|^2, \, d_4(x)/|x|^4\}$, where $d_2(x)=\max \{ |x_4^2+x_5^2|, \, |x_0|^2, \, |x_1|^2, \, |x_2|^2, \, |x_3|^2 \}$ and $d_4(x)=\max \{ |x_1^4+x_2^4+x_3^4|, \, |x_0|^4, \, |x_4|^4, \, |x_5|^4 \}$ for \texttt{CICY1} and \texttt{CICY2}. The vertical axis shows $\displaystyle \max_{v \in T_xM : |v|_{FS}=1} |v|_{h}$, where $|\cdot|_{FS}$ denotes the length of a tangent vector in the ambient Fubini-Study metric and $|\cdot|_h$ denotes the length of a vector in the learned approximate Calabi-Yau metric. For \texttt{Quintic} and \texttt{CICY2}, which are small perturbations of a singular variety, the data shows a long neck forming. \texttt{CICY1} is a large perturbation of a singular variety and the neck formation cannot be observed.}
    \label{figure:neck-forming}
\end{figure}

The data is presented in \cref{figure:neck-forming}.

\subsection{Approximate harmonic $1$-forms}
\label{section:results-1-forms}

\paragraph{Comparison of four example varieties.}
We studied the four varieties \texttt{Fermat}, \texttt{Quintic}, \texttt{CICY1}, and \texttt{CICY2}.
\begin{itemize}
    \item 
    The real locus of \texttt{Fermat} has first Betti number equal to zero, so its only harmonic $1$-form is constant zero.
    Consequently, when trying to numerically solve the harmonic $1$-form equation, we expect a poor quality solution.

    \item 
    The real locus of \texttt{Quintic} does admit a non-zero harmonic $1$-form, but it must have at least one zero by \cref{proposition:must-have-zeros-on-RP3-with-handle}.

    \item 
    The real locus of \texttt{CICY2} is diffeomorphic to $S^1 \times S^2$, to it admits a harmonic $1$-form.
    Because \texttt{CICY2} is close to a singular limit whose real locus admits a $\text{U}(1)$-symmetry, we speculated that \texttt{CICY2} may admit a nowhere vanishing harmonic $1$-form.

    \item 
    The real locus of \texttt{CICY1} is also diffeomorphic to $S^1 \times S^2$, but because \texttt{CICY1} is not close to a singular Calabi-Yau manifold, we have no reason to expect that the (unique up to scaling) harmonic $1$-form is nowhere vanishing.
\end{itemize}

\begin{figure}
    \centering
    \includegraphics[width=9cm]{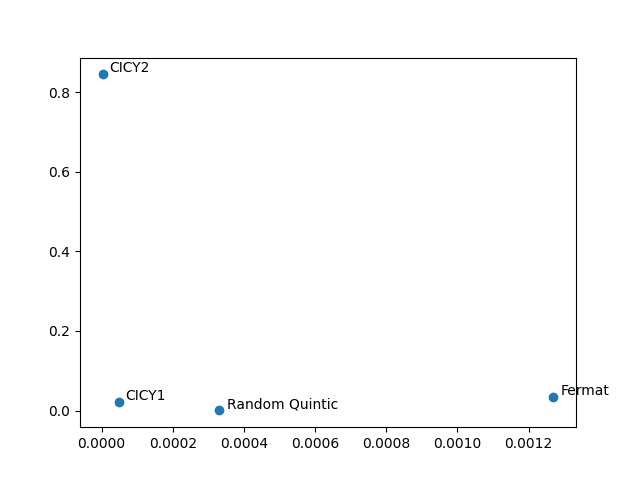}
    \caption{Plot of the data from \cref{tab:harmonic_loss}. $x$-axis shows the approximation error of the approximate $1$-form $\omega$, $y$-axis shows the minimum value of $|\omega |$ over the manifold. The data confirms our expectation and neatly fits into a pattern of large/small approximation error and having/not having a zero.}
    \label{fig:harmonic-loss-four-dots-plot}
\end{figure}

We see these expectations confirmed by the experimental data presented in \cref{tab:harmonic_loss} and \cref{fig:harmonic-loss-four-dots-plot}.

\begin{table}[h]
    \centering
    \begin{tabular}{|c|c|c|}
        \hline
        Manifold & Harmonic loss & Normalized min norm\\
        \hline
        \texttt{Fermat} & $1.27 \times 10^{-3}$ & $3.30 \times 10^{-2}$\\
        \texttt{Quintic} & $3.31 \times 10^{-4}$ & $2.88 \times 10^{-4}$\\
        \texttt{CICY1} & $4.75 \times 10^{-5}$ & $2.25 \times 10^{-2}$\\
        \texttt{CICY2} & $2.42 \times 10^{-6}$ & $8.45 \times 10^{-1}$\\
        \hline
    \end{tabular}
    \caption{Results of approximately computing a harmonic $1$-form on four example varieties. \emph{Harmonic loss} measures the quality of the approximate solution. On varieties whose real locus has zero first Betti number, we expect this quality to be low. \emph{Normalised min norm} is the minimum norm in a large number of sampled points of the approximate $1$-form scaled to have unit $L^2$-norm. On some varieties, we expect the harmonic $1$-form to have no zeros, on others we have shown it must have zeros. This is explained in \cref{section:results-1-forms}.}
    \label{tab:harmonic_loss}
\end{table}

\paragraph{\texttt{CICY2}.}
The data in \cref{fig:harmonic-loss-four-dots-plot,tab:harmonic_loss} suggest that there exists a nowhere vanishing harmonic $1$-form on the real locus $L=S^1 \times S^2$ of \texttt{CICY2}.
We speculated in \cref{subsection:intersections-of-quadrics-and-quartics} that this is due to a $\text{U}(1)$-action on the nearby singular Calabi-Yau variety $X^*$, and that the harmonic $1$-form should be approximately constant in the $S^1$-direction and approximately zero in the $S^2$ direction.

To see if this is the case, we plot the $1$-form restricted to:
\begin{itemize}
    \item 
    $L \cap \{(x_0:x_1:x_2:x_3:x_4:x_5):x_0=1,x_1=0,x_2=0 \} \subset \R^3$, which is approximately $Z(C_{\text{aff}}) \times Z(x_3^4) \cong S^1 \times \{\pm 1\} \subset \R^3$,
    where we expect the $1$-form to be approximately constant in the $S^1$-direction here;

    \item 
    $L \cap \{(x_0:x_1:x_2:x_3:x_4:x_5):x_0=1,x_1=0,x_4=0 \} \subset \R^3$, which is approximately $Z(x_2^4+x_3^4) \times Z(x_5^2) \cong S^1 \times \{\pm 1\}$, where we expect the $1$-form to be approximately zero.
\end{itemize}

\begin{figure}[htbp]
    \centering
    \begin{minipage}{0.45 \textwidth}
        \centering
        \includegraphics[width=8cm]{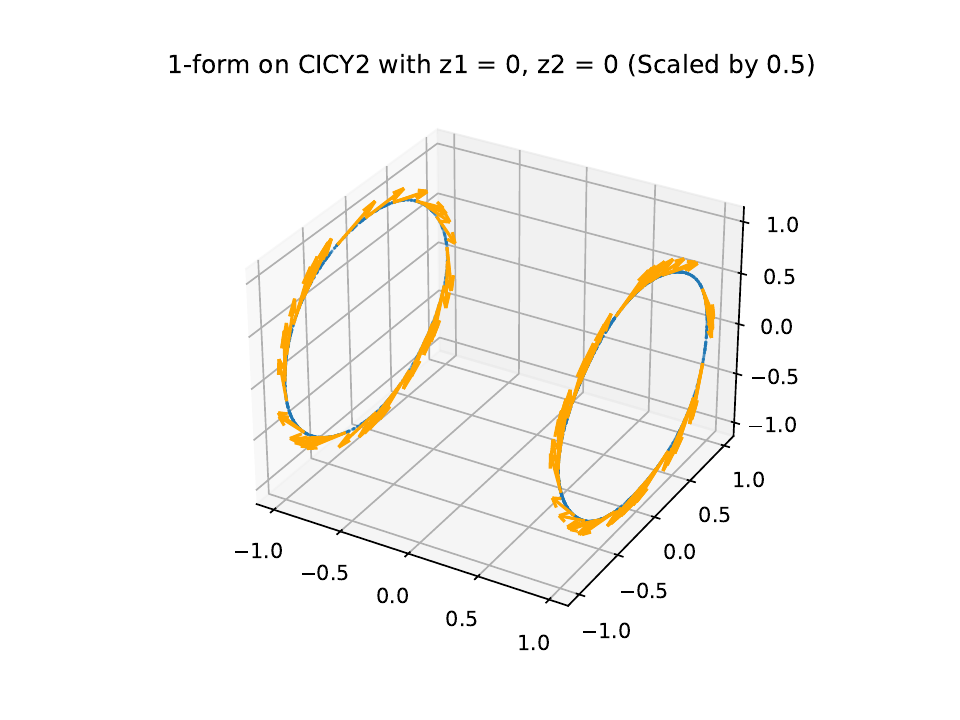}
    \end{minipage}
    \begin{minipage}{0.45 \textwidth}
        \centering
        \includegraphics[width=8cm]{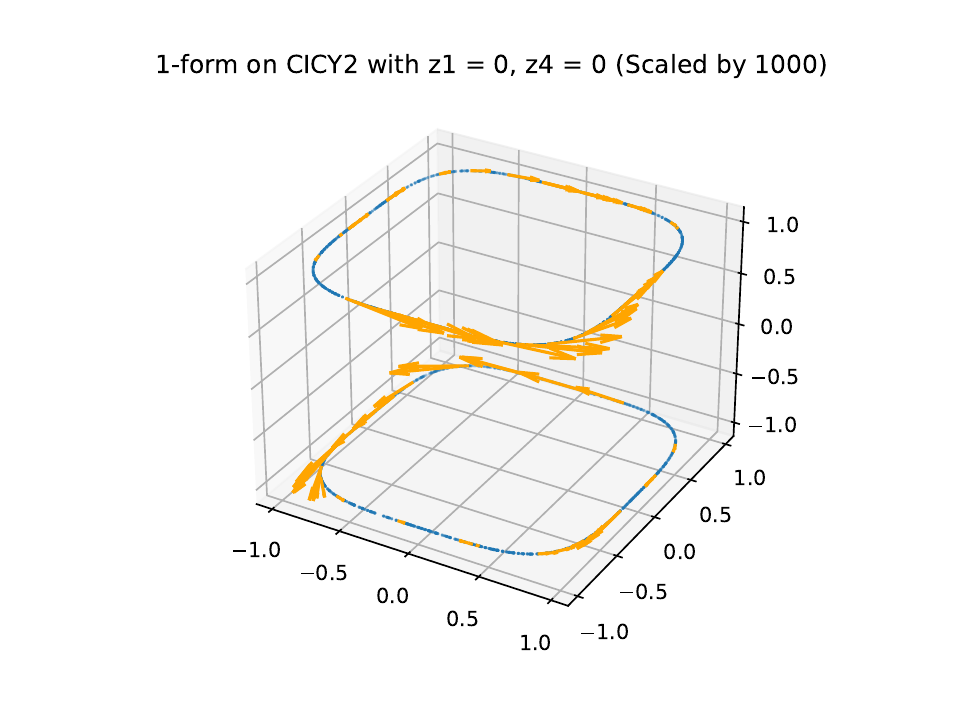}
    \end{minipage}
    \caption{The approximate harmonic $1$-form restricted to different one-dimensional submanifolds of \texttt{CICY2}.
    For better visualisation, arrow lengths have been rescaled by a factor of $0.5$ in the left picture and by a factor of $1000$ in the right picture.
    The plots confirm the expectation that the approximate harmonic $1$-form is approximately constant in the $S^1$-direction.}
    \label{figure:cicy2-circle-plots}
\end{figure}

We see this expectation confirmed in \cref{figure:cicy2-circle-plots}.
We therefore make the following conjecture:

\begin{conjecture}
    Let $X_\epsilon$ be a $1$-parameter deformation of the singular variety $X^*$ from \cref{subsection:intersections-of-quadrics-and-quartics}, and denote by $g_{\text{CY}}^\epsilon$ the uniquely determined Calabi-Yau metric in Kähler class $\omega_{\text{FS}} \mid _{X_\epsilon}$.
    Then, we have
    \[
        g_{\text{CY}}^\epsilon \mid_L
        \rightarrow
        g_{\text{product}}
        \text{ in }
        C^\infty
        \text{ as }
        \epsilon \rightarrow 0,
    \]
    where $g_{\text{product}}$ denotes a product metric on $S^1 \times S^2$.
\end{conjecture}

\begin{figure}[htbp]
    \centering
    \begin{minipage}{0.45 \textwidth}
        \centering
        \includegraphics[width=8cm]{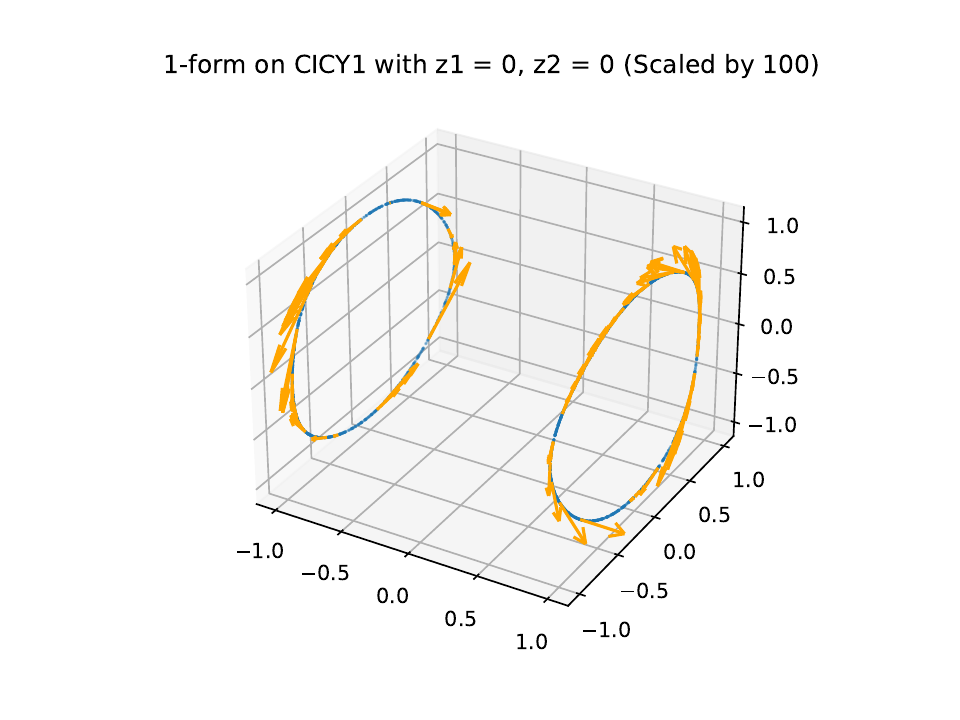}
    \end{minipage}
    \begin{minipage}{0.45 \textwidth}
        \centering
        \includegraphics[width=8cm]{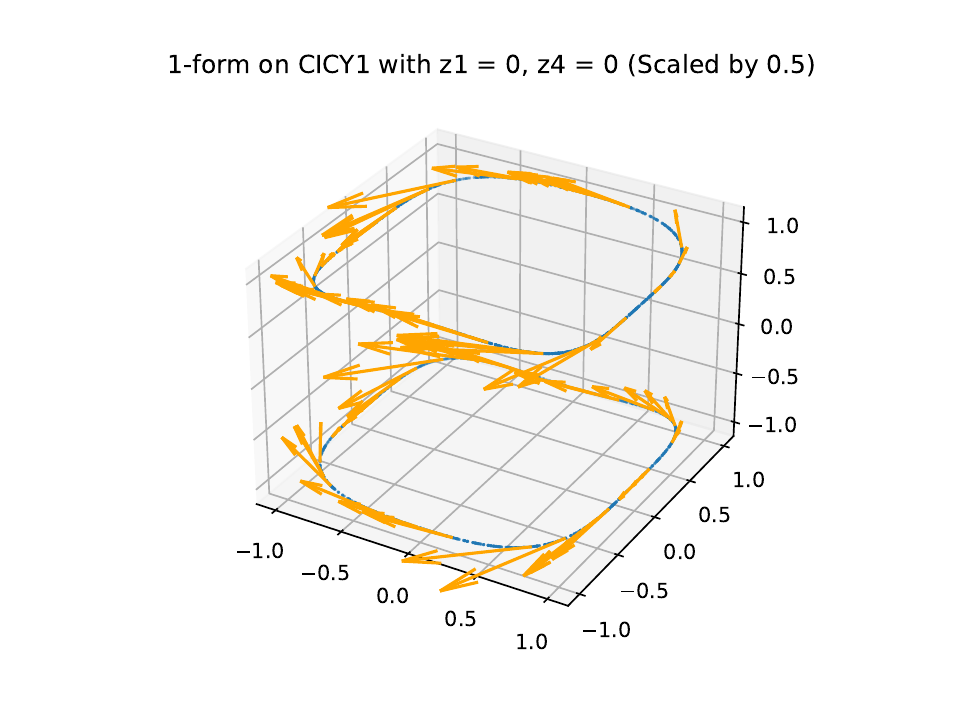}
    \end{minipage}

    \begin{minipage}{0.45 \textwidth}
        \centering
        \includegraphics[width=8cm]{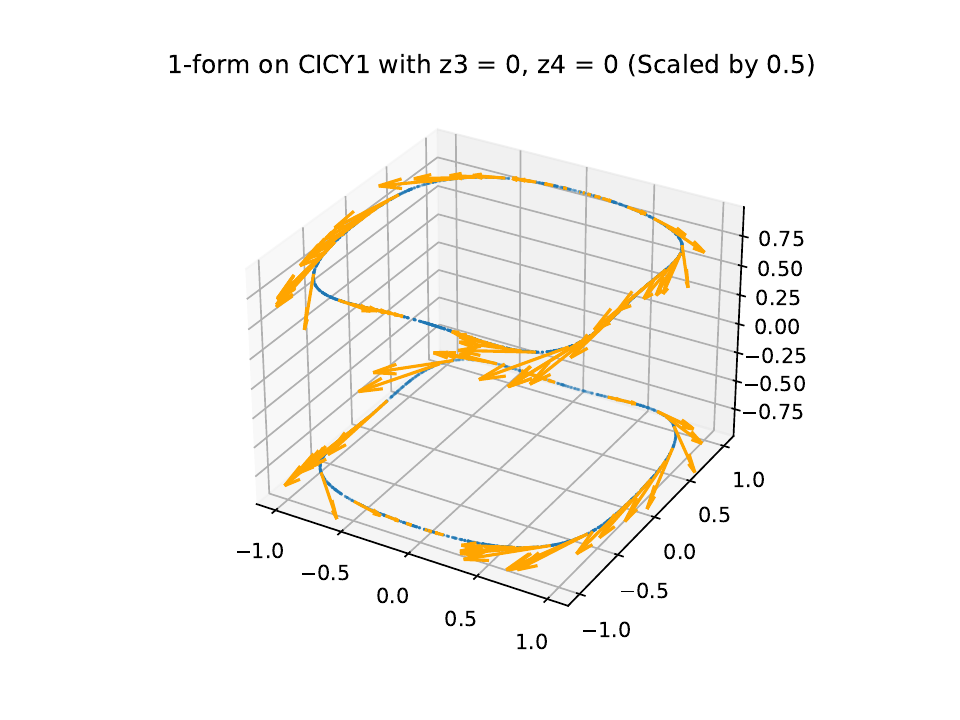}
    \end{minipage}
    \begin{minipage}{0.45 \textwidth}
        \centering
        \includegraphics[width=8cm]{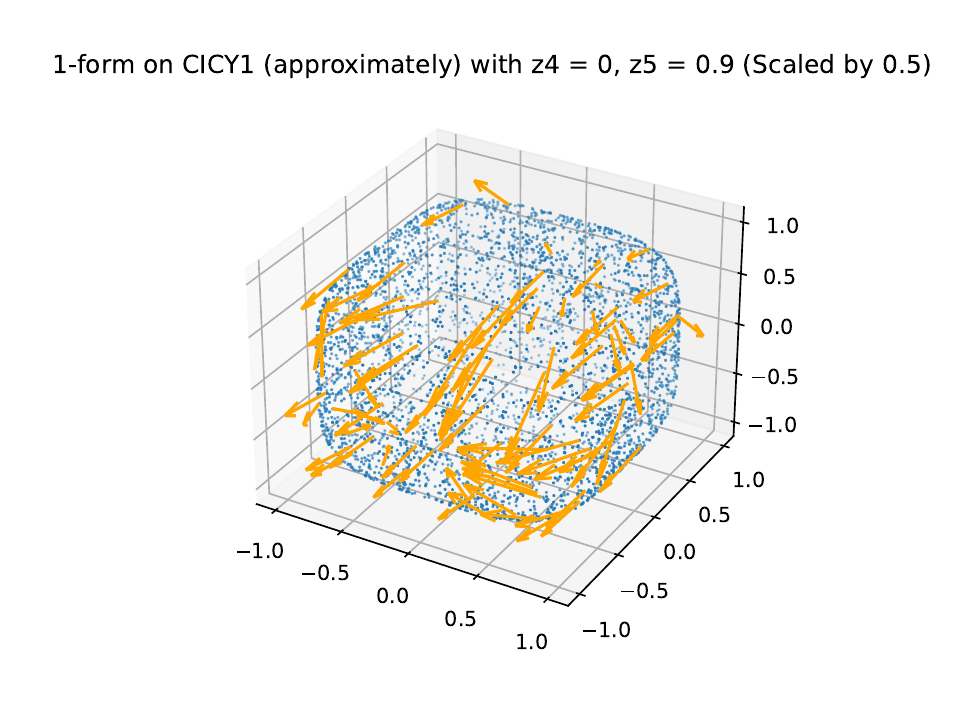}
    \end{minipage}
    \caption{First three pictures show the approximate harmonic $1$-form restricted to different one-dimensional submanifolds of \texttt{CICY1}.
    Fourth picture shows the $1$-form on the variety $-1+z_1^4+z_2^4+z_3^4=0$.
    (Note that this variety is \emph{not} contained in \texttt{CICY1}, but only nearby.)
    For better visualisation, arrow lengths have been rescaled as follows: $100$ (top-left), $0.5$ (remaining three pictures).
    The plots show that the $1$-form is approximately zero in the $S^1$-direction, and on $S^2$ appears to have a single sink and a single source.}
    \label{figure:cicy1-circle-sphere-plots}
\end{figure}

\paragraph{\texttt{CICY1}.}
The data in \cref{fig:harmonic-loss-four-dots-plot,tab:harmonic_loss} suggest that there exists a harmonic $1$-form on the real locus $L \cong S^1 \times S^2$ of \texttt{CICY1}.
Inspection of the approximate solution shows that it is approximately zero in the $S^1$-direction, and approximately a $1$-form flowing from a single source to a single sink in the $z_2$-direction on $S^2$.
This is visualised in \cref{figure:cicy1-circle-sphere-plots}.
We note that such a $1$-form on $S^2$ would be approximately the differential of an eigenfunction of the Laplacian with minimum eigenvalue, so it minimises the functional $(\|{\d \omega}_{L^2}+\|{\d^* \omega}_{L^2})/\|{\omega}_{L^2}$ among forms on $S^2$.
This result is somewhat unexpected, but we found that the neural network learned this solution, even when initialised to be the solution from \texttt{CICY2}, which is far from non-zero in the $S^1$-direction.

\appendix

\section{Implementation Details}
\label{section:implementation-details}
Table \ref{tab:metric-training} shows the hyperparameters for the training of Calabi-Yau metrics. 
The points on the manifolds are sampled using the algorithms described in \ref{subsection:sampling-algorithms},
with the ratio of training points to testing points being approximately 10:1.
To avoid numerical instability caused by dividing by a small number,
the points are further distributed into patches $U_{i, j}$, where $i = \arg\max_{i} \left| z^i \right|$
and $j = \arg\max_{j \neq i} \left| \partial_j f\right|$. 
(For \texttt{CICY1} and \texttt{CICY2} this becomes $U_{i, \{j, k\}}$ with 
$j, k = \arg\max_{j \neq k \neq i} \left| \det \partial (f_1, f_2) / \partial (z^j, z^k) \right|$).
We then use a 3-hidden-layer (equivalently degree 8 polynomials) bihomogeneous neural network
(Eq. \ref{eq:bihomo-net}) to represent the Kähler potential 
and compute the Kähler metric on $\mathbb{CP}^4$.
After restricting the metric to the manifold, we train the neural network with the MAPE loss:
\begin{align}\label{eq:MAPE-loss}
\int_Y
\left| \frac{\omega^3}{\Omega \wedge \bar{\Omega}}-1 \right|
\operatorname{vol},
\end{align}
using Adam algorithm followed by 100 epochs of L-BFGS.
Since \texttt{Quintic} has more terms and is less symmetric, 
we performed a more comprehensive hyperparameter tuning and trained it for more epochs. 

\begin{table}[h]
    \centering
    \begin{tabular}{|c|c|c|c|}
        \hline
        Manifold & Train pts & Layers & \# Epochs\\
        \hline
        \texttt{Fermat} & 500000 & (64, 256, 1024, 1) & 4000\\
        \texttt{Quintic} & 100000 & (256, 512, 512, 1) & 11430 \\
        \texttt{CICY1} & 84183 & (64, 256, 1024, 1)  & 4000\\
        \texttt{CICY2} & 84393 & (64, 256, 1024, 1)  & 4000\\
        \hline
    \end{tabular}
    \caption{The number of training points, the widths of the network and the number of epochs for the training of the Ricci-flat metrics on each manifold.} 
    \label{tab:metric-training}
\end{table}

\begin{figure}[htbp]
    \centering
    \begin{minipage}{0.24 \textwidth}
        \centering
        \includegraphics[width=3.5cm]{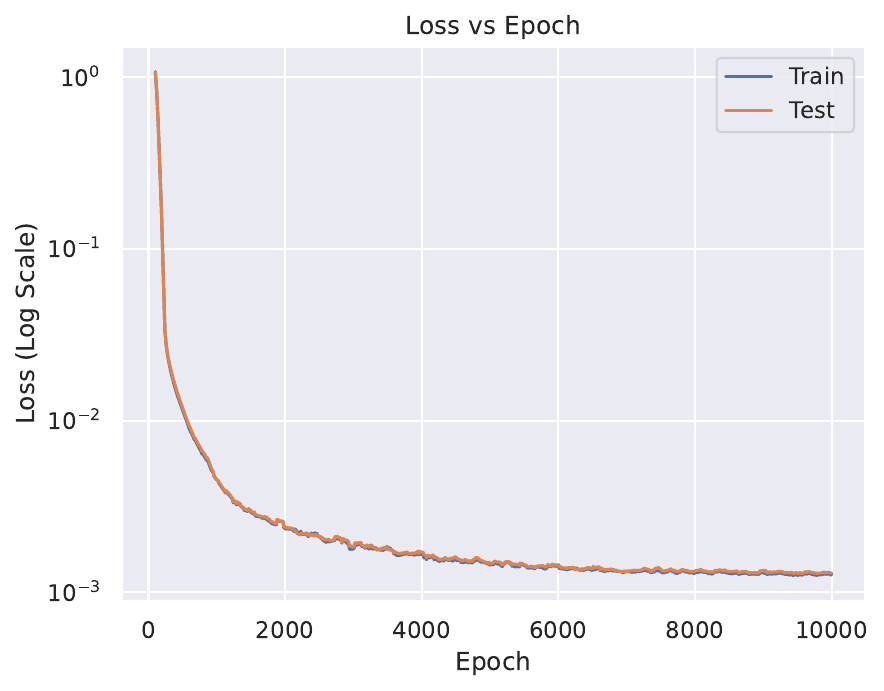}
    \end{minipage}
    \begin{minipage}{0.24 \textwidth}
        \centering
        \includegraphics[width=3.5cm]{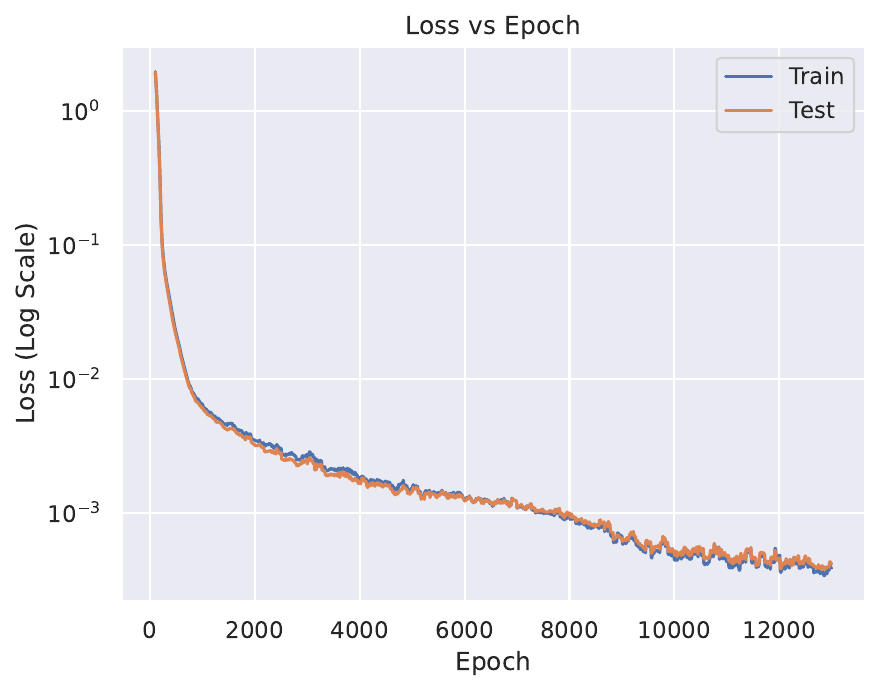}
    \end{minipage}
    \begin{minipage}{0.24 \textwidth}
        \centering
        \includegraphics[width=3.5cm]{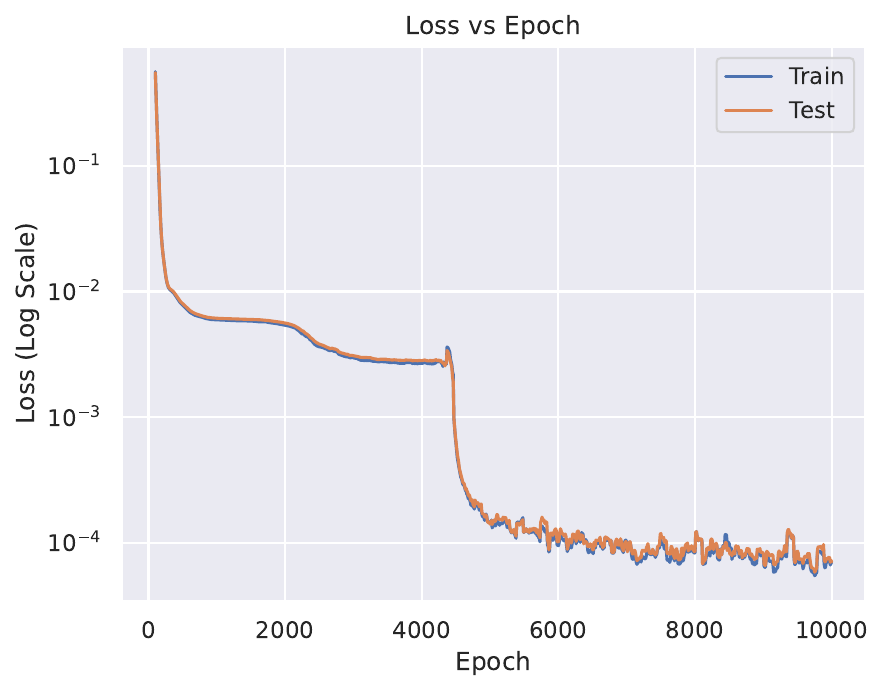}
    \end{minipage}
    \begin{minipage}{0.24 \textwidth}
        \centering
        \includegraphics[width=3.5cm]{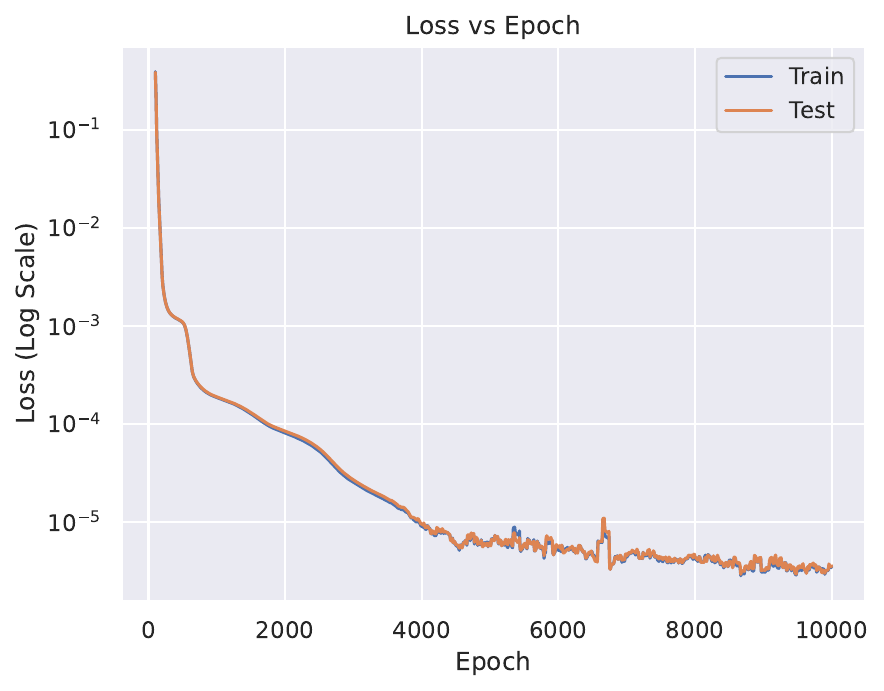}
    \end{minipage}
    \caption{The train and test loss curve for the harmonic 1-forms. From left to right: \texttt{Fermat}, \texttt{Quintic}, \texttt{CICY1}, \texttt{CICY2}.}
    \label{figure:loss-curve}
\end{figure}

\begin{figure}[htbp]
    \centering
    \begin{minipage}{0.24 \textwidth}
        \centering
        \includegraphics[width=3.5cm]{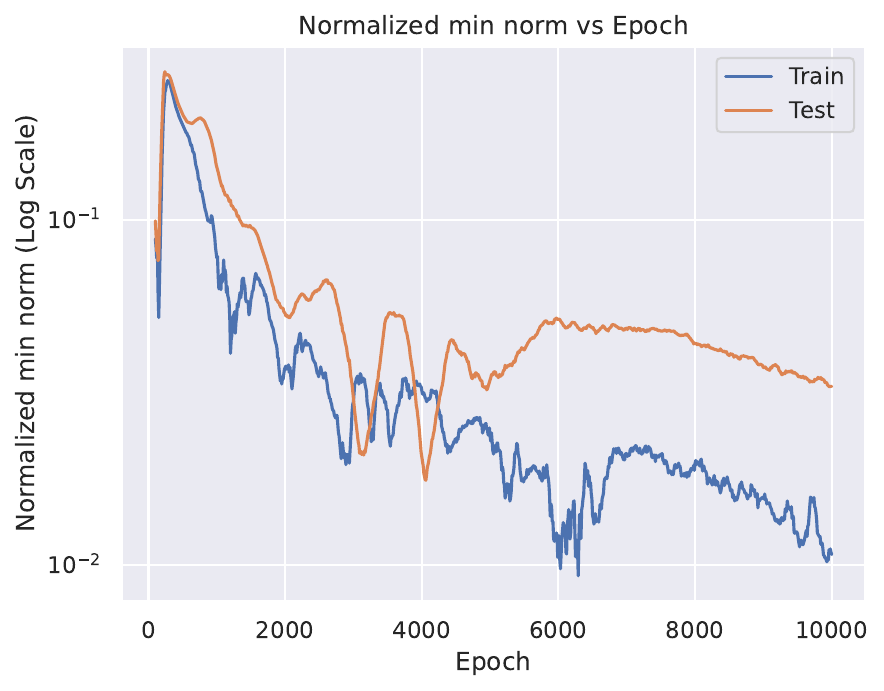}
    \end{minipage}
    \begin{minipage}{0.24 \textwidth}
        \centering
        \includegraphics[width=3.5cm]{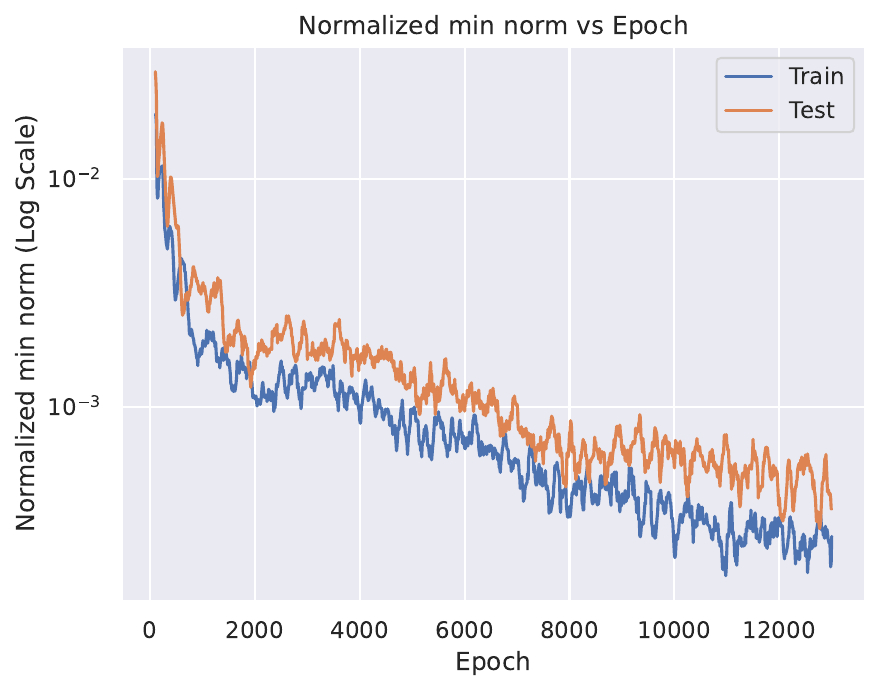}
    \end{minipage}
    \begin{minipage}{0.24 \textwidth}
        \centering
        \includegraphics[width=3.5cm]{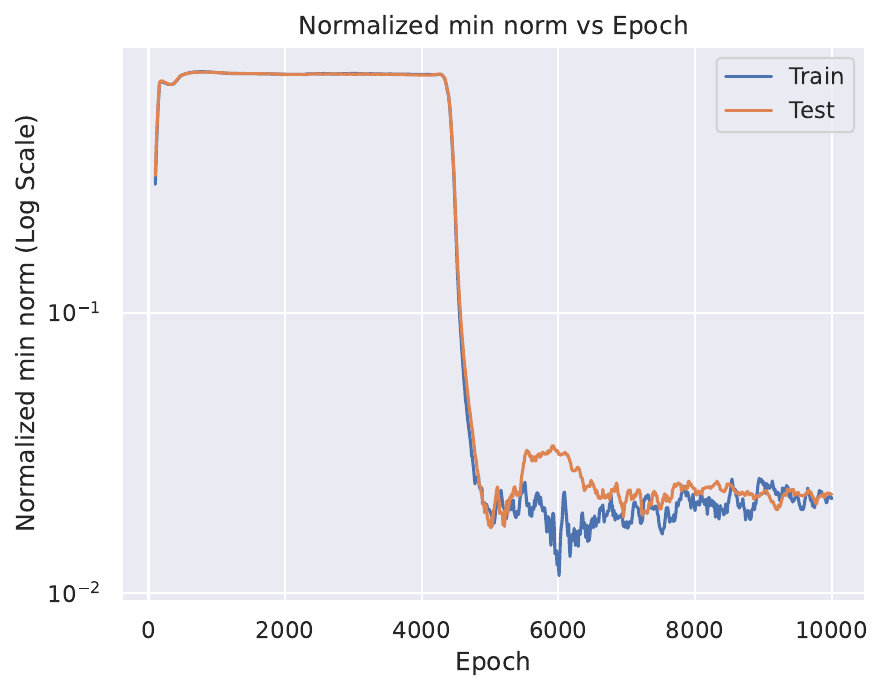}
    \end{minipage}
    \begin{minipage}{0.24 \textwidth}
        \centering
        \includegraphics[width=3.5cm]{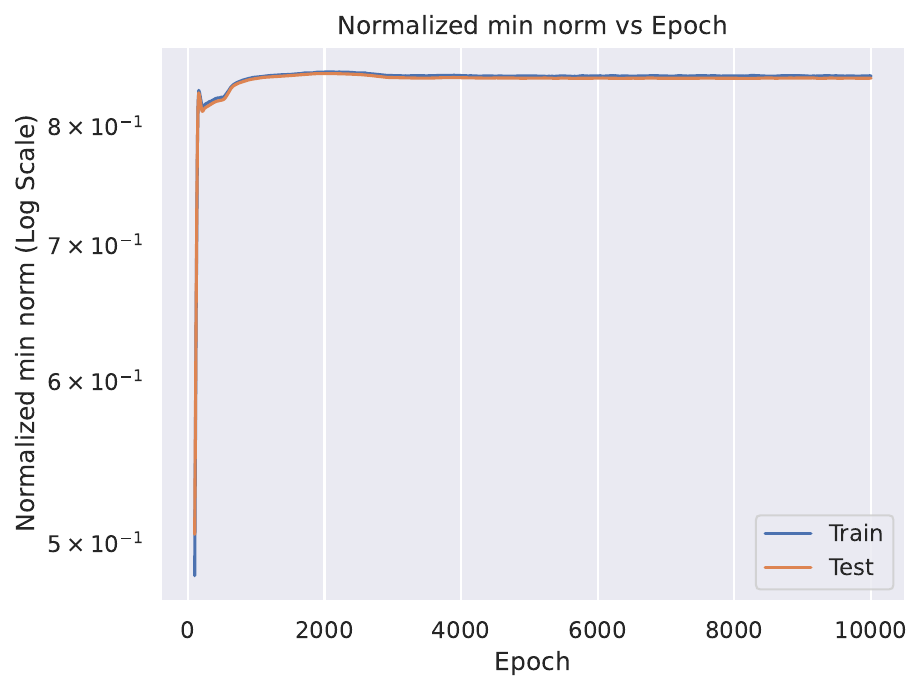}
    \end{minipage}
    \caption{The minimal norm curve for the harmonic 1-forms, normalized so that the average norms are scaled to one on each manifold (Eq. \ref{eq:norm-constraint}). From left to right: \texttt{Fermat}, \texttt{Quintic}, \texttt{CICY1}, \texttt{CICY2}.}
    \label{figure:min-norm-curve}
\end{figure}

Table \ref{tab:one-form-training} shows the hyperparameters for the training of harmonic 1-forms. The points on the special Lagrangian submanifolds are sampled and distributed into patches using the same algorithms for Calabi-Yaus, and the 1-forms are constructed and trained using the method described in \ref{section:numerical-forms}. The discrepancies between train sets and test sets are insignificant as shown in Fig. \ref{figure:loss-curve} and Fig. \ref{figure:min-norm-curve} 
\begin{table}[h]
    \centering
    \begin{tabular}{|c|c|c|c|}
        \hline
        Manifold & Train pts & Layers & \# Epochs\\
        \hline
        \texttt{Fermat} & 125596 & (128, 256, 1024, 10) & 10000\\
        \texttt{Quintic} & 139829 & (128, 256, 1024, 10) & 13000\\
        \texttt{CICY1} & 42053 & (128, 256, 1024, 15)  & 10000\\
        \texttt{CICY2} & 43366 & (128, 256, 1024, 15)  & 10000\\
        \hline
    \end{tabular}
    \caption{The number of training points, the widths of the network and the number of epochs for the training of the harmonic 1-forms on each manifold. \texttt{CICY1} and \texttt{CICY2} have wider output layers since they are embedded in a higher dimensional space.} 
    \label{tab:one-form-training}
\end{table}

\pagebreak
\bibliographystyle{apalike}
\bibliography{library}

\begin{thebibliography}{}

\bibitem[Anderson et~al., 2010]{Anderson2010}
Anderson, L.~B., Braun, V., Karp, R.~L., and Ovrut, B.~A. (2010).
\newblock Numerical {H}ermitian {Y}ang-{M}ills connections and vector bundle
  stability in heterotic theories.
\newblock {\em Journal of High Energy Physics}, 2010(6):1--45.

\bibitem[Anderson et~al., 2021]{Anderson2021}
Anderson, L.~B., Gerdes, M., Gray, J., Krippendorf, S., Raghuram, N., and
  Ruehle, F. (2021).
\newblock Moduli-dependent {C}alabi-{Y}au and ${SU}(3)$-structure metrics from
  machine learning.
\newblock {\em Journal of High Energy Physics}, 2021(5):1--45.

\bibitem[Ashmore et~al., 2020]{Ashmore2020}
Ashmore, A., He, Y.~H., and Ovrut, B.~A. (2020).
\newblock Machine {Learning} {Calabi}–{Yau} {Metrics}.
\newblock {\em Fortschritte der Physik}, 68(9).

\bibitem[Ashmore and Ruehle, 2021]{Ashmore2021}
Ashmore, A. and Ruehle, F. (2021).
\newblock Moduli-dependent {KK} towers and the swampland distance conjecture on
  the quintic {C}alabi-{Y}au manifold.
\newblock {\em Physical Review D}, 103(10):106028.

\bibitem[Berglund et~al., 2022]{Berglund2022}
Berglund, P., Butbaia, G., H{\"u}bsch, T., Jejjala, V., Pe{\~n}a, D.~M.,
  Mishra, C., and Tan, J. (2022).
\newblock Machine learned {C}alabi-{Y}au metrics and curvature.
\newblock {\em arXiv preprint arXiv:2211.09801}.

\bibitem[Bochner, 1946]{Bochner1946}
Bochner, S. (1946).
\newblock Vector fields and {R}icci curvature.
\newblock {\em Bull. Amer. Math. Soc.}, 52:776--797.

\bibitem[Boykov and Kolmogorov, 2003]{Boykov2003}
Boykov and Kolmogorov (2003).
\newblock Computing geodesics and minimal surfaces via graph cuts.
\newblock In {\em Proceedings Ninth IEEE international conference on computer
  vision}, pages 26--33. IEEE.

\bibitem[Braun et~al., 2008]{Braun2008}
Braun, V., Brelidze, T., Douglas, M.~R., and Ovrut, B.~A. (2008).
\newblock {C}alabi-{Y}au metrics for quotients and complete intersections.
\newblock {\em Journal of High Energy Physics}, 2008(05):080.

\bibitem[Breiding and Marigliano, 2020]{Breiding2020}
Breiding, P. and Marigliano, O. (2020).
\newblock Random points on an algebraic manifold.
\newblock {\em SIAM Journal on Mathematics of Data Science}, 2(3):683--704.

\bibitem[Buchbinder et~al., 2014]{Buchbinder2014}
Buchbinder, E.~I., Constantin, A., and Lukas, A. (2014).
\newblock The moduli space of heterotic line bundle models: a case study for
  the tetra-quadric.
\newblock {\em Journal of High Energy Physics}, 2014(3):1--46.

\bibitem[Bunch and Donaldson, 2008]{Bunch2008}
Bunch, R.~S. and Donaldson, S.~K. (2008).
\newblock Numerical approximations to extremal metrics on toric surfaces.
\newblock In {\em Handbook of geometric analysis. {N}o. 1}, volume~7 of {\em
  Adv. Lect. Math. (ALM)}, pages 1--28. Int. Press, Somerville, MA.

\bibitem[Butbaia et~al., 2024]{Butbaia2024}
Butbaia, G., Pe{\~n}a, D.~M., Tan, J., Berglund, P., H{\"u}bsch, T., Jejjala,
  V., and Mishra, C. (2024).
\newblock Physical {Y}ukawa couplings in heterotic string compactifications.
\newblock {\em arXiv preprint arXiv:2401.15078}.

\bibitem[Candelas et~al., 1987]{Candelas1987}
Candelas, P., Dale, A.~M., L\"{u}tken, C.~A., and Schimmrigk, R.~C. (1987).
\newblock Complete intersection {C}alabi-{Y}au manifolds.
\newblock In {\em Frontiers of high energy physics ({L}ondon, 1986)}, pages
  88--134. Hilger, Bristol.

\bibitem[Constantin et~al., 2024]{Constantin2024}
Constantin, A., Fraser-Taliente, C.~S., Harvey, T.~R., Lukas, A., and Ovrut, B.
  (2024).
\newblock Computation of quark masses from string theory.
\newblock {\em arXiv preprint arXiv:2402.01615}.

\bibitem[Di~Rocco et~al., 2022]{DiRocco2022}
Di~Rocco, S., Eklund, D., and G{\"a}fvert, O. (2022).
\newblock Sampling and homology via bottlenecks.
\newblock {\em Mathematics of Computation}, 91(338):2969--2995.

\bibitem[Donaldson, 2009]{Donaldson2009}
Donaldson, S.~K. (2009).
\newblock Some numerical results in complex differential geometry.
\newblock {\em Pure Appl. Math. Q.}, 5(2):571--618.

\bibitem[Doran et~al., 2008]{Doran2008}
Doran, C., Headrick, M., Herzog, C.~P., Kantor, J., and Wiseman, T. (2008).
\newblock Numerical {K}\"{a}hler-{E}instein metric on the third del {P}ezzo.
\newblock {\em Comm. Math. Phys.}, 282(2):357--393.

\bibitem[Douglas et~al., 2022]{Douglas2022}
Douglas, M., Lakshminarasimhan, S., and Qi, Y. (2022).
\newblock Numerical {C}alabi-{Y}au metrics from holomorphic networks.
\newblock In Bruna, J., Hesthaven, J., and Zdeborova, L., editors, {\em
  Proceedings of the 2nd Mathematical and Scientific Machine Learning
  Conference}, volume 145 of {\em Proceedings of Machine Learning Research},
  pages 223--252. PMLR.

\bibitem[Douglas, 2022]{Douglas2022b}
Douglas, M.~R. (2022).
\newblock Holomorphic feedforward networks.
\newblock {\em Pure Appl. Math. Q.}, 18(1):251--268.

\bibitem[Douglas et~al., 2008]{Douglas2008}
Douglas, M.~R., Karp, R.~L., Lukic, S., and Reinbacher, R. (2008).
\newblock Numerical {C}alabi-{Y}au metrics.
\newblock {\em J. Math. Phys.}, 49(3):032302, 19.

\bibitem[Gerdes and Krippendorf, 2023]{gerdes_cyjax_2023}
Gerdes, M. and Krippendorf, S. (2023).
\newblock {CYJAX}: {A} package for {Calabi}-{Yau} metrics with {JAX}.
\newblock {\em Machine Learning: Science and Technology}, 4(2):025031.
\newblock arXiv:2211.12520 [hep-th].

\bibitem[G\'{o}mez-Serrano, 2019]{Gomez-Serrano2019}
G\'{o}mez-Serrano, J. (2019).
\newblock Computer-assisted proofs in {PDE}: a survey.
\newblock {\em SeMA J.}, 76(3):459--484.

\bibitem[Gäfvert, 2022]{Gafvert2022}
Gäfvert, O. (2022).
\newblock Sampling and homology via bottlenecks.
\newblock \url{
  https://www.JuliaHomotopyContinuation.org/examples/sampling_bottlenecks/ }.
\newblock Accessed: January 21, 2022.

\bibitem[Headrick and Nassar, 2013]{Headrick2013}
Headrick, M. and Nassar, A. (2013).
\newblock Energy functionals for {C}alabi-{Y}au metrics.
\newblock In {\em Journal of Physics: Conference Series}, volume 462, page
  012019. IOP Publishing.

\bibitem[Headrick and Wiseman, 2005]{Headrick2005}
Headrick, M. and Wiseman, T. (2005).
\newblock Numerical {R}icci-flat metrics on {K3}.
\newblock {\em Classical and Quantum Gravity}, 22(23):4931.

\bibitem[Jaco, 1980]{Jaco1980}
Jaco, W. (1980).
\newblock {\em Lectures on three-manifold topology}, volume~43 of {\em CBMS
  Regional Conference Series in Mathematics}.
\newblock American Mathematical Society, Providence, RI.

\bibitem[Jejjala et~al., 2022]{jejjala_neural_2022}
Jejjala, V., Peña, D. K.~M., and Mishra, C. (2022).
\newblock Neural network approximations for {Calabi}-{Yau} metrics.
\newblock {\em Journal of High Energy Physics}, 2022(8):105.

\bibitem[Joyce and Karigiannis, 2021]{Joyce2021}
Joyce, D. and Karigiannis, S. (2021).
\newblock A new construction of compact torsion-free {$\rm G_2$}-manifolds by
  gluing families of {E}guchi-{H}anson spaces.
\newblock {\em J. Differential Geom.}, 117(2):255--343.

\bibitem[Joyce, 1996]{Joyce1996}
Joyce, D.~D. (1996).
\newblock Compact {R}iemannian {$7$}-manifolds with holonomy {$G_2$}. {I},
  {II}.
\newblock {\em J. Differential Geom.}, 43(2):291--328, 329--375.

\bibitem[Joyce, 2000]{Joyce2000}
Joyce, D.~D. (2000).
\newblock {\em Compact manifolds with special holonomy}.
\newblock Oxford Mathematical Monographs. Oxford University Press, Oxford.

\bibitem[Krasnov, 2006]{Krasnov2006}
Krasnov, V.~A. (2006).
\newblock Rigid isotopy classification of real three-dimensional cubics.
\newblock {\em Izv. Ross. Akad. Nauk Ser. Mat.}, 70(4):91--134.

\bibitem[Krasnov, 2009]{Krasnov2009}
Krasnov, V.~A. (2009).
\newblock Topological classification of real three-dimensional cubics.
\newblock {\em Mathematical Notes}, 85(5):841--847.

\bibitem[Larfors et~al., 2021]{larfors2021learning}
Larfors, M., Lukas, A., Ruehle, F., and Schneider, R. (2021).
\newblock Learning size and shape of calabi-yau spaces.

\bibitem[Larfors et~al., 2022]{Larfors2022}
Larfors, M., Lukas, A., Ruehle, F., and Schneider, R. (2022).
\newblock Numerical metrics for complete intersection and {K}reuzer--{S}karke
  {C}alabi--{Y}au manifolds.
\newblock {\em Machine Learning: Science and Technology}, 3(3):035014.

\bibitem[Legland et~al., 2007]{Legland2007}
Legland, D., Ki{\^e}u, K., and Devaux, M.-F. (2007).
\newblock Computation of {M}inkowski measures on 2d and 3d binary images.
\newblock {\em Image Analysis and Stereology}, 26(2):83--92.

\bibitem[Lehmann and Legland, 2012]{Lehmann2012}
Lehmann, G. and Legland, D. (2012).
\newblock Efficient n-dimensional surface estimation using {C}rofton formula
  and run-length encoding.
\newblock {\em Efficient N-Dimensional surface estimation using Crofton formula
  and run-length encoding, Kitware INC (2012)}.

\bibitem[Li et~al., 2003]{Li2003}
Li, X., Wang, W., Martin, R.~R., and Bowyer, A. (2003).
\newblock Using low-discrepancy sequences and the {C}rofton formula to compute
  surface areas of geometric models.
\newblock {\em Computer-Aided Design}, 35(9):771--782.

\bibitem[Li, 2023]{Li2023}
Li, Y. (2023).
\newblock Intermediate complex structure limit for {C}alabi-{Y}au metrics.

\bibitem[Milnor, 1965]{Milnor1965}
Milnor, J.~W. (1965).
\newblock {\em Topology from the differentiable viewpoint}.
\newblock University Press of Virginia, Charlottesville, VA.
\newblock Based on notes by David W. Weaver.

\bibitem[Morozov, 2022]{Morozov2022}
Morozov, E. (2022).
\newblock Surfaces containing two isotropic circles through each point (v1).

\bibitem[Qi, 2022]{QiGithub2024}
Qi, Y. (2022).
\newblock {G}it{H}ub - yidiq7/{M}{L}{G}eometry: {M}achine learning
  {C}alabi-{Y}au metrics --- github.com.
\newblock \url{https://github.com/yidiq7/MLGeometry/tree/master}.
\newblock [Accessed 01-05-2024].

\bibitem[Reiterer and Trubowitz, 2019]{Reiterer2019}
Reiterer, M. and Trubowitz, E. (2019).
\newblock Choptuik's critical spacetime exists.
\newblock {\em Comm. Math. Phys.}, 368(1):143--186.

\bibitem[Seyyedali, 2009]{Seyyedali2009}
Seyyedali, R. (2009).
\newblock Numerical algorithm for finding balanced metrics on vector bundles.
\newblock {\em Asian J. Math.}, 13(3):311--321.

\bibitem[Shiffman and Zelditch, 1999]{shiffman_distribution_1999}
Shiffman, B. and Zelditch, S. (1999).
\newblock Distribution of zeros of random and quantum chaotic sections of
  positive line bundles.
\newblock {\em Communications in Mathematical Physics}, 200(3):661--683.
\newblock arXiv: math/9803052.

\bibitem[Sun, 2022]{Sun2022}
Sun, S. (2022).
\newblock Geometry of {C}alabi-{Y}au metrics.
\newblock {\em Notices Amer. Math. Soc.}, 69(4):546--556.

\bibitem[Sun and Zhang, 2019]{Sun2019}
Sun, S. and Zhang, R. (2019).
\newblock Complex structure degenerations and collapsing of {C}alabi-{Y}au
  metrics.
\newblock {\em arXiv preprint arXiv:1906.03368}.

\bibitem[Tischler, 1970]{Tischler1970}
Tischler, D. (1970).
\newblock On fibering certain foliated manifolds over {$S^1$}.
\newblock {\em Topology}, 9:153--154.

\bibitem[Wu, 2017]{Wu2017}
Wu, H.-H. (2017).
\newblock {\em The {B}ochner technique in differential geometry}, volume~6 of
  {\em CTM. Classical Topics in Mathematics}.
\newblock Higher Education Press, Beijing, new edition.
\newblock Expanded version of [MR0714349].

\bibitem[Yau, 1977]{Yau1977}
Yau, S.~T. (1977).
\newblock Calabi's conjecture and some new results in algebraic geometry.
\newblock {\em Proc. Nat. Acad. Sci. U.S.A.}, 74(5):1798--1799.

\bibitem[Yau, 1978]{Yau1978}
Yau, S.~T. (1978).
\newblock On the {R}icci curvature of a compact {K}\"{a}hler manifold and the
  complex {M}onge-{A}mp\`ere equation. {I}.
\newblock {\em Comm. Pure Appl. Math.}, 31(3):339--411.

\bibitem[Yau, 2009]{Yau2009}
Yau, S.-T. (2009).
\newblock A survey of {C}alabi-{Y}au manifolds.
\newblock In {\em Surveys in differential geometry. {V}ol. {XIII}. {G}eometry,
  analysis, and algebraic geometry: forty years of the {J}ournal of
  {D}ifferential {G}eometry}, volume~13 of {\em Surv. Differ. Geom.}, pages
  277--318. Int. Press, Somerville, MA.

\end{thebibliography}

\end{document}